\documentclass[a4paper]{article}
\usepackage{amsmath}

\title{Some binomial identities related to Calabi-Yau differential equations}
\author{Gert Almkvist.}
\begin{document}

\maketitle
In [1] we found some 300 new differential equations of Calabi-Yau type. Most
of them were found by using Maple's \textquotedblright
Zeilberger\textquotedblright\ on the coefficients $A_{n}$ of the the
analytic solution $y_{0}=\sum_{n=0}^{\infty }A_{n}z^{n}$ . Sometimes rather
different sums of products of binomial coefficients gave the same
recursionformula and also the same $A_{n}.$These cases are listed here with
the same numbering as in [1]. Even more often did it happen for double sums
(here we did not use \textquotedblright Zeilberger\textquotedblright ). All
identities containing only simple sums are proved by the Zeilberger
machinery. $A_{n}^{\prime }$ denotes coefficients of fifth order equations.
At the end of the paper we also list some formulas for $A_{n}$ obtained as
the constant term of the n-th power of a Laurent polynomial coming from a
reflexive polytope in \textbf{R}$^{4}.$ These formulas were found in
collaboration with Victor Batyrev, Max Kreuzer and Duco van Straten.

We also list some ''ghost formulas'' for $"A_n"$ where ''Zeilberger'' gives
the same recursion formula as for $A_n$ but $"A_n"$ is certainly not equal
to $A_n.$ Here is an example of this phenomen:

Let 
\[
"A_n"=\left\{ \frac{(3n)!}{n!^3}\right\}
^b\sum_k(-1)^{n+k(b+c)}k^{-c}(k+\frac n2)(k-n)^b\binom{k+2n}{3n}^b\binom{n+k}%
n^{-c} 
\]
\[
\]

which usually is divergent. Then $y_0=\sum_{n=0}^\infty "A_n"z^n$
''satisfies'' the following differential equations ($\theta =z\frac d{dz}$ ) 
\[
\]

$%
\begin{array}{rrr}
b & c & \text{differential equation} \\ 
1 & 4 & \theta ^2-z(11\theta ^2+11\theta +3)+z^2(\theta +1)^2 \\ 
1 & 5 & \theta ^3+z(2\theta +1)(17\theta ^2+17\theta +5)+z^2(\theta +1)^3 \\ 
1 & 6 & \#198 \\ 
1 & 7 & \text{diff. eq. of }\zeta (5).\text{and }\zeta (3)\text{, see [2]
(10.4)} \\ 
&  &  \\ 
2 & 3 & \theta ^2-3z^2(3\theta +2)(3\theta +4) \\ 
2 & 4 & \theta ^3-z(2\theta +1)(13\theta ^2+13\theta +4)-3z^2(\theta
+1)(3\theta +2)(3\theta +4) \\ 
2 & 5 & \#195 \\ 
2 & 6 & \#32%
\end{array}
$

In [2] it is explained that the correct formula for the coefficients of the
solutions of the differential equations above are 
\[
A_n=\left\{ \frac{(3n)!}{n!^3}\right\} ^b\sum_{k=0}^n\binom nk^c\binom{3n}{%
n+k}^{-b}\left\{ 1+k(-cH_k+cH_{n-k}+bH_{n+k}-bH_{2n-k})\right\} 
\]
where 
\[
H_k=\sum_{j=1}^n\frac 1j 
\]
and $H_k=0$ if $k\leq 0.$ In [2] most formulas of type 
\[
"A_n"=\sum_k(n-2k)C(n,k) 
\]
where $C(n,n-k)=C(n,k)$ are converted to harmonic sums and are proved. 
\[
\]

First we note the following very simple identity 
\[
\sum_k\binom nk^2\binom{3k}{2n}=\sum_k\binom nk^2\binom{2k}k 
\]
which deserves a bijective proof. Not so simple is that these sums also are
equal to 
\[
3^{1-n}\binom{2n}n\sum_k(-1)^k\frac{n-2k}{2n-3k}\binom nk^2\binom{3n-3k}{2n}%
\binom{2n}{3k}^{-1} 
\]

Here are three other identities that do not belong to the differential
equations in the table. 
\[
\sum_k(-1)^{n+k}\binom nk\binom{n+k}n^2=\sum_k\binom nk^2\binom{n+k}n 
\]
\[
\sum_k(-1)^{n+k}\binom nk\binom{n+k}n\binom{2n-2k}{n-k}=\frac{(4p)!}{%
(2p)!p!^2}\text{ if }n=3p\text{, }0\text{ otherwise} 
\]
\[
\sum_k(-1)^{n+k}\binom nk\binom{n+k}n\binom{2n-k}n=\frac{(3p)!}{p!^3}\text{
if }n=2p,\text{ }0\text{ otherwise } 
\]

\textbf{3.} 
\[
A_n=\binom{2n}n^4 
\]
\[
A_n=\binom{2n}n\sum_k\binom nk^2\binom{n+k}k\binom{2n+k}n 
\]
\[
A_n=\binom{2n}n\sum_{i,j}\binom ni^2\binom nj\binom{i+j}j\binom{2n}j 
\]
\[
\]

\textbf{15.} 
\[
A_n=\frac{(3n)!}{n!^3}\sum_k\binom nk^3 
\]
\[
A_n=\binom{3n}n\sum_k\binom nk\binom{n+k}n\binom{2n-2k}n\binom{2n}{n+k} 
\]
\[
A_n=2^{-n}\binom{2n}n\sum_k\binom nk\binom{n+k}n\binom{2k}n\binom{3n}{n+k}%
\binom{2n}{2k}^{-1} 
\]
\[
A_n=\sum_k\binom nk^2\binom{n+k}k\binom{2n-k}n\binom{3n}{n+k} 
\]
\[
A_n=\sum_k\binom nk\binom{n+k}n\binom{2n-k}n\binom{3n}{n+k}\binom{2k}n 
\]
\[
A_n=\sum_{i,j}\binom ni\binom nj\binom{2n-i}n\binom{2i}n\binom{2n}i\binom{2n}%
j 
\]
\[
\]

\textbf{16. } 
\[
A_n=\binom{2n}n\sum_{i+j+k+l=n}\left\{ \frac{n!}{i!j!k!l!}\right\} ^2 
\]
\[
A_n=\binom{2n}n\sum_k\binom nk^2\binom{2k}k\binom{2n-2k}{n-k} 
\]
\[
A_n=\sum_k\binom{2n}{2k}\binom{2k}k^2\binom{2n-2k}{n-k}^2 
\]
\[
A_n=4^{-n}\binom{2n}n\sum_k\binom nk^{-2}\binom{2k}k^3\binom{2n-2k}{n-k}^3 
\]
\[
A_n=\binom{2n}n^2\sum_k(-)^{n+k}\binom nk^2\binom{2k}k\binom k{n-k}\binom{2n%
}{2k}^{-1} 
\]
\[
A_n=\binom{2n}n\sum_k(-1)^{n+k}\binom nk\binom{2k}k\binom{2n-2k}{n-k}\binom{%
2k}n 
\]
\[
A_n=\sum_k\binom nk\binom{2n-k}n\binom{2k}k\binom{2n-2k}{n-k}\binom{2n}k 
\]
\[
A_n=\sum_k\binom nk\binom{n+k}n\binom{2k}k\binom{2n-2k}{n-k}\binom{2n}{n-k} 
\]
\[
A_n=\binom{2n}n\sum_{i,j}\binom ni^2\binom nj^2\binom{2j}j 
\]
\[
A_n=\binom{2n}n\sum_{i,j}\binom ni\binom ij^3\binom{2n-i}n 
\]
\[
A_n=(-1)^n\binom{2n}n\sum_{i,j}(-8)^{n-i}\binom ni\binom ij^3\binom{n+i}n 
\]
\[
A_n=(-1)^n\binom{2n}n\sum_{i,j}(-8)^{n-i}\binom ni\binom ij^3\binom{2i}n 
\]
\[
A_n=(-1)^n\binom{2n}n\sum_{i,j}(-8)^{n-i}\binom ni\binom ij^2\binom{n+i}n%
\binom{2j}i 
\]
\[
A_n=\binom{2n}n\sum_{i,j}(-1)^{n+i+j}4^{n-i-j}\binom ni\binom ij\binom{2i}i%
\binom{2j}j\binom{n+i+j}n 
\]
\[
A_n=\binom{2n}n\sum_{i,j}(-1)^{n+i+j}4^{n-i-j}\binom ni\binom ij\binom{2i}i%
\binom{2j}j\binom{n+i+j}{n+j} 
\]
\[
A_n=\binom{2n}n\sum_{i,j}(-1)^{n+i+j}4^{n-i-j}\binom ni\binom ij\binom{2i}i%
\binom{2j}j\binom{n+i+j}{n+j} 
\]
\[
A_n=\binom{2n}n\sum_{i,j}(-1)^{n+i+j}4^{n-i-j}\binom ni\binom ij\binom{2i}i%
\binom{2j}j\binom{2n-i-j}n 
\]
\[
\]

\textbf{18. } 
\[
A_n=\binom{2n}n\sum_k\binom nk^4 
\]
\[
A_n=\binom{2n}n\sum_k\binom{2k}k\binom{2n-2k}{n-k}\binom k{n-k}\binom{2n-k}k 
\]
\[
A_n=\binom{2n}n\sum_k\binom nk\binom{2k}k\binom{2n-2k}n\binom{n+k}{n-k} 
\]
\[
A_n=\binom{2n}n\sum_k\binom nk^2\binom{2k}n\binom{2n-2k}n 
\]
\[
A_n=\binom{2n}n\sum_k\binom nk\binom{2k}k\binom{2k}n\binom k{n-k} 
\]
\[
A_n=\binom{2n}n^2\sum_k(-1)^k\binom nk^2\binom{2k}k\binom{n+k}{n-k}\binom{2n%
}{2k}^{-1} 
\]
\[
A_n=\binom{2n}n^2\sum_k(-1)^k\binom nk^3\binom{2n-k}n\binom{2n}{2k}^{-1} 
\]
\[
A_n=\binom{2n}n^2\sum_k(-1)^{n+k}\binom nk^2\binom{2k}k\binom{2n-2k}{n-k}%
\binom{2n}k^{-1} 
\]
\[
A_n=\binom{2n}n\sum_k(-1)^{n+k}\binom nk^2\binom{n+k}n\binom{2n-k}n\binom{2n}%
k\binom{2n}{2k}^{-1} 
\]
\[
A_n=\binom{2n}n\sum_k(-1)^k\binom nk^2\binom{2n-k}n^2\binom{2n}k\binom{2n}{2k%
}^{-1} 
\]
\[
A_n=\binom{2n}n\sum_{i,j}\binom ni^2\binom nj\binom ij\binom{2j}n 
\]
\[
A_n=\binom{2n}n\sum_{i,j}\binom ni^2\binom nj\binom ij\binom i{n-j} 
\]
\[
A_n=\binom{2n}n\sum_{i,j}\binom ni^2\binom nj\binom ij\binom{n-j}i 
\]
\[
A_n=\binom{2n}n\sum_{i,j}\binom ni\binom nj\binom{i+j}j\binom{2i}n\binom
n{i+j} 
\]
\[
A_n=\binom{2n}n\sum_{i,j}\binom ni\binom nj\binom ij\binom{n+i-j}n\binom
j{i-j} 
\]
\[
A_n=\binom{2n}n\sum_{i,j}(-1)^j4^{n-j}\binom ni^2\binom nj\binom{i+j}j\binom{%
2j}j 
\]
\[
A_n=\binom{2n}n\sum_{i,j}(-1)^j4^{n-j}\binom ni^2\binom nj\binom{i+j}n\binom{%
2j}j 
\]
\[
A_n=\binom{2n}n\sum_{i,j}(-1)^{i+j+n}\binom ni\binom nj\binom{2n-i}n\binom{%
2n-j}n\binom{i+j}n 
\]
\[
A_n=\binom{2n}n\sum_{i,j}(-1)^{i+j+n}\binom ni\binom nj\binom{2n-i}n\binom{%
2n-j}n\binom{i+j}j 
\]
\[
A_n=\binom{2n}n\sum_k(-1)^{n+k}\binom nk\binom{n+k}n\binom{2k}k\binom{2n-2k}{%
n-k} 
\]
\[
"A_n"=\binom{2n}n\sum_kk^{-4}\binom{n+k}n^{-4} 
\]
\[
"A_n"=\binom{2n}n\sum_kk^{-3}(k-n)^2\binom kn^3\binom{n+k}n^{-1} 
\]
\[
\]

\textbf{19.} 
\[
A_n=\sum_k\binom nk^3\binom{n+k}n\binom{2n-k}n 
\]
\[
A_n=\sum_k\binom nk^2\binom{2n-k}n\binom{2k}k\binom{n+k}{n-k} 
\]
\[
A_n=\sum_{i,j}\binom nj^2\binom nj^2\binom ij\binom{n+j}n 
\]
\[
A_n=\sum_{i,j}\binom ni^2\binom nj\binom ij\binom{n+i}n\binom{2j}n 
\]
\[
A_n=\sum_{i,j}\binom ni^2\binom nj\binom ij\binom{i+j}n\binom{n+i}{n-j} 
\]
\[
A_n=\sum_{i,j}\binom ni^2\binom nj\binom{i+j}j\binom{n-j}i\binom{n+i}{n-j} 
\]
\[
A_n=\sum_{i,j}\binom ni^2\binom nj\binom{n+j}n\binom{2j}n\binom{n+i}{n+j} 
\]
\[
A_n=\sum_{i,j}\binom ni^2\binom nj\binom ij\binom i{n-j}\binom{n+i}n 
\]
\[
A_n=\sum_{i,j}\binom ni^2\binom nj^2\binom{n+i-j}{n-j}\binom{2j}n 
\]
\[
A_n=\sum_{i,j}\binom ni\binom nj\binom{2n-i}n\binom{n+j}n\binom{2j}n\binom
n{i-j} 
\]
\[
A_n=\sum_{i,j}\binom ni^2\binom nj\binom{n+i}n\binom{i+j}j\binom{n+i}{n+j} 
\]
\[
A_n=\sum_{i,j}\binom ni^2\binom nj\binom{n+j}n\binom{2j}n\binom n{i+j} 
\]
\[
A_n=\sum_{i,j}\binom ni^2\binom nj\binom{n+j}n\binom{2j}n\binom{2n-i}{n+j} 
\]
\[
A_n=\sum_{i,j}\binom ni\binom nj\binom{n+i}n\binom{n+j}n\binom{2j}n\binom
n{i+j} 
\]
\[
A_n=\sum_{i,j}\binom ni\binom nj\binom{2n-j}n\binom{2i}n\binom{i+j}j\binom
n{i+j} 
\]
\[
A_n=\sum_{i,j}\binom ni^2\binom nj^2\binom{n+i-j}i\binom{2j}n 
\]
\[
A_n=\sum_{i,j}\binom ni^2\binom{n+j}n^2\binom{i+j}n\binom{3n+1}{n-2j} 
\]
\[
A_n=\sum_{i,j}\binom ni^2\binom{n+j}n^2\binom{i+j}n\binom{2n-i-j}n 
\]
\[
A_n=\sum_{i,j}\binom ni^2\binom{n+j}n^2\binom{i+j}n\binom{2n-i-j}{n-j} 
\]
\[
A_n=\sum_{i,j}\binom ni^2\binom nj\binom ij\binom{i+j}j\binom{2j}n 
\]
\[
A_n=\sum_{i,j}\binom ni^2\binom nj\binom ij\binom{i+j}n\binom{n+i-j}{n-j} 
\]
\[
\]

\textbf{21.} 
\[
A_n=\sum_k\binom nk^3\binom{2k}k\binom{2n-2k}{n-k} 
\]
\[
A_n=\sum_{i,j}\binom ni^2\binom nj\binom ij^2\binom{2j}j 
\]
\[
"A"_n=\sum_k(n-2k)\binom nk^5\binom{2k}n\binom{2n-2k}n 
\]
\[
\]

\textbf{22.} 
\[
A_n=\sum_k\binom nk^5 
\]
\[
A_n=\sum_{i,j}\binom ni^2\binom nj^2\binom ij\binom i{n-j} 
\]
\[
A_n=\sum_{i,j}\binom ni^2\binom nj^2\binom{2i}n\binom{2n-i-j}n 
\]
\[
A_n=\sum_{i,j}\binom ni^2\binom nj^2\binom{n+i-j}n\binom{2j}n 
\]
\[
A_n=\sum_{i,j}(-1)^{i+j}\binom{2n-i}n^3\binom{n+j}n^2\binom{2n+1}{i-j} 
\]
\[
A_n=\sum_{i,j}\binom ni^2\binom nj\binom{2n-i}n\binom ij\binom i{n-j} 
\]
\[
A_n=\sum_{i,j}\binom ni^2\binom nj\binom ij\binom{2j}n\binom{n+i-j}{n-j} 
\]
\[
A_n=\sum_{i,j}\binom ni^2\binom nj\binom ij\binom{2j}n\binom{2n-i}n 
\]
\[
A_n=\sum_{i,j}\binom ni\binom nj\binom{n+j}n\binom{i+j}j\binom{2i}n\binom
n{i+j} 
\]
\[
"A_n"=\sum_k(-1)^kk^{-5}\binom{n+k}n^{-5} 
\]
\[
\]

\textbf{24.} 
\[
A_n=\frac{(3n)!}{n!^3}\sum_k\binom nk^2\binom{n+k}n 
\]
\[
A_n=\sum_k(-1)^{n+k}\binom{n+k}n^3\binom{2n-k}n\binom{3n}{n+k} 
\]
\[
"A"_n=\binom{2n}n^2\sum_k(-1)^kk^{-3}\binom nk^{-1}\binom{n+k}n^{-2} 
\]
\[
"A_n"=\left\{ \frac{(3n)!}{n!^3}\right\} ^2\sum_k(-1)^{n+k}(n-k)k^{-2}\binom
nk\binom{2n+k}{3n}\binom{n+k}n^{-3} 
\]
\[
\]

\textbf{25.} 
\[
A_n=\binom{2n}n^2\sum_k\binom nk^2\binom{n+k}n 
\]
\[
A_n=\binom{2n}n^2\sum_k\binom nk\binom{2k}k\binom{n+k}{n-k} 
\]
\[
A_n=\binom{2n}n\sum_k\binom nk\binom{n+k}n\binom{2n-k}n\binom{2n}k 
\]
\[
A_n=\binom{2n}n\sum_k\binom nk\binom{n+k}n\binom{2n-k}n\binom{2n}{n-k} 
\]
\[
A_n=\binom{2n}n\sum_k\binom nk^2\binom{2k}n\binom{n+2k}n 
\]
\[
A_n=\binom{2n}n^3\sum_k(-1)^{n+k}\binom nk^2\binom{2n-k}n\binom{2n}k^{-1} 
\]
\[
A_n=\binom{2n}n^4\sum_k(-1)^k\binom nk^3\binom{2n}k^{-2} 
\]
\[
A_n=\binom{2n}n^2\sum_k(-1)^k\binom nk\binom{n+k}n^2 
\]
\[
A_n=\sum_{i,j}\binom ni^2\binom nj^2\binom{i+j}j\binom{3n-i}n 
\]
\[
A_n=\sum_{i.j}\binom ni^2\binom nj\binom{n+i}n\binom{n+i-j}n\binom{2n}{n-j} 
\]
\[
A_n=\sum_{i,j}\binom ni\binom nj\binom{2n-i}n\binom{2n+i-j}{2n}\binom{2n}i%
\binom{2n}j 
\]
\[
A_n=\binom{2n}n\sum_{i,j}\binom ni^2\binom nj\binom{i+j}j\binom{n+i}{n-j} 
\]
\[
A_n=\binom{2n}n\sum_{i,j}\binom ni^2\binom nj\binom{2n+j}{2n}\binom{n+i}{n+j}
\]
\[
A_n=\binom{2n}n\sum_{i,j}\binom ni^2\binom nj\binom{n+i-j}{n-j}\binom{n+i}{%
n-j} 
\]
\[
A_n=\binom{2n}n\sum_{i,j}\binom ni^2\binom nj\binom ij\binom{2n-i-j}i 
\]
\[
A_n=\binom{2n}n\sum_{i,j}\binom ni^2\binom nj\binom{n+i-j}{n-j}\binom{n+i}{%
n+j} 
\]
\[
A_n=\binom{2n}n\sum_{i,j}\binom ni^2\binom nj\binom{n+i-j}{n-j}\binom{n+i}j 
\]
\[
A_n=\binom{2n}n^2\sum_{i,j}(-1)^{i+j}\binom{2n-i}n^2\binom{2n-j}n\binom{2n+1%
}{i-j} 
\]
\[
A_n=\binom{2n}n^2\sum_{i,j}(-1)^{i+j}\binom{2n-i}n^2\binom{n+j}n\binom{2n+1}{%
i-j} 
\]
\[
A_n=\binom{2n}n^2\sum_{i,j}(-1)^{i+j}\binom ni\binom nj\binom{n+j}n^2\binom{%
n+i-j}n 
\]
\[
A_n=\binom{2n}n^2\sum_{i,j}(-1)^{i+j}\binom ni\binom nj\binom{n+j}n^2\binom{%
2n+i+j}n 
\]
\[
A_n=\sum_{i,j}\binom ni\binom nj\binom{i+j}n^3\binom{2n}{i+j} 
\]
\[
A_n=\sum_{i,j}\binom ni\binom nj\binom{i+j}n^2\binom{3n-i-j}n\binom{2n}{i+j} 
\]
\[
A_n=\sum_{i,j}(-1)^{i+j}\binom{2n-i}n^2\binom{n+j}n^2\binom{2n+j}n\binom{3n+1%
}{i-j} 
\]
\[
A_n=\sum_{i,j}(-1)^{i+j}\binom{2n-i}n^2\binom{n+j}n^2\binom{3n-i}n\binom{3n+1%
}{i-j} 
\]
\[
A_n=\sum_{i,j}(-1)^{i+j}\binom{2n-i}n^2\binom{n+j}n^2\binom{3n-i}n\binom{3n+1%
}{i-j} 
\]
\[
A_n=\binom{2n}n\sum_{i,j}(-1)^{i+j}\binom{2n-i}n^2\binom{n+j}n\binom{2n+j}n%
\binom{3n+1}{n-j} 
\]
\[
A_n=\binom{2n}n\sum_{i,j}(-1)^j\binom ni^2\binom{n+j}n\binom{n+i}n\binom{3n+1%
}{n-j} 
\]
\[
A_n=\binom{2n}n\sum_{i,j}(-1)^j\binom ni^2\binom{n+j}n\binom{2n-i}n\binom{%
3n+1}{n-j} 
\]
\[
A_n=\binom{2n}n\sum_{i,j}(-1)^{i+j}\binom ni\binom nj\binom{i+j}n^3 
\]
\[
A_n=\binom{2n}n^2\sum_{i,j}(-1)^{i+j}\binom ni\binom nj\binom{n+j}n\binom
n{2j-i} 
\]
\[
A_n=\binom{2n}n^2\sum_{i,j}(-1)^{i+j}\binom ni\binom nj\binom{n+2j}{2n}%
\binom n{2j-i} 
\]
\[
A_n=\binom{2n}n^2\sum_{i,j}(-1)^{i+j}\binom ni\binom nj\binom{n+j}n^2\binom{%
n+i+j}n 
\]
\[
A_n=\binom{2n}n^2\sum_{i,j}\binom ni\binom ij\binom{2j}j\binom{2j}{i-j} 
\]
\[
A_n=\binom{2n}n\sum_{i,j}\binom ni^2\binom nj\binom{2i}n\binom{2i}{n-j} 
\]
\[
A_n=\binom{2n}n^2\sum_{i,j}\binom ni^2\binom nj\binom{2i}{n+j} 
\]
\[
A_n=\binom{2n}n\sum_{i,j}\binom ni^2\binom nj\binom ij\binom{2n+i-j}i 
\]
\[
A_n=\binom{2n}n\sum_{i,j}\binom ni\binom nj\binom{2n-i}n\binom{n+i}n\binom
n{i+j} 
\]
\[
A_n=\sum_{i,j}\binom ni\binom nj\binom{n+i}n\binom{n+j}n\binom{2n+j}n\binom
n{i+j} 
\]
\[
A_n=\binom{2n}n^2\sum_{i,j}\binom ni\binom nj\binom{n-i}j\binom n{i+j} 
\]
\[
A_n=\binom{2n}n^2\sum_{i,j}\binom ni\binom nj\binom{i+j}j\binom n{i+j} 
\]
\[
A_n=\binom{2n}n^2\sum_{i,j}(-1)^{i+j}\binom{n+i}n\binom{n+j}n\binom{i+j}%
j\binom n{i+j} 
\]
\[
A_n=\binom{2n}n^2\sum_{i,j}(-1)^{i+j}\binom{2n-i}n\binom{n+j}n\binom{i+j}%
j\binom n{i+j} 
\]
\[
A_n=\sum_{i,j}\binom ni\binom nj\binom{2n-i}n\binom{n+j}n\binom{2n+j}n\binom
n{i-j} 
\]
\[
A_n=\sum_{i,j}\binom ni\binom nj\binom{2n-i}n\binom{n+j}n\binom{3n-i}n\binom
n{i-j} 
\]
\[
A_n=\binom{2n}n\sum_{i,j}\binom ni\binom nj\binom{n+j}n^2\binom n{i-j} 
\]
\[
"A_n"=\binom{2n}n^2\sum_kk^{-3}\binom nk^{-2}\binom{n+k}n^{-1} 
\]
\[
"A_n"=(2n+1)^{-2}\sum_k\binom{n+k}n^3\binom{2n+k+1}k^{-2} 
\]
\[
"A_n"=(2n+1)^{-1}\binom{2n}n\sum_k\binom{n+k}n^3\binom{2n+k+1}k^{-1} 
\]
\[
\]

\textbf{26.} 
\[
A_n=\binom{2n}n\sum_k\binom nk^2\binom{n+k}n\binom{2k}n 
\]
\[
A_n=\sum_k\binom nk\binom{n+k}n\binom{2n-k}n\binom{2n}k\binom{2k}n 
\]
\[
A_n=\binom{2n}n\sum_k\binom nk^2\binom{2k}n\binom{2n-k}n 
\]
\[
A_n=\binom{2n}n\sum_k\binom nk^2\binom{n+k}k\binom{2n-k}n 
\]
\[
A_n=\binom{2n}n\sum_k\binom nk\binom{n+k}n\binom{2n-2k}{n-k}\binom{2n-k}k 
\]
\[
A_n=\binom{2n}n\sum_k\binom nk\binom{n+k}n\binom{2k}k\binom k{n-k} 
\]
\[
A_n=\binom{2n}n\sum_k\binom nk\binom{2k}k\binom{2k}n\binom{n+k}{n-k} 
\]
\[
A_n=\binom{2n}n\sum_k\binom nk\binom{2k}k\binom{2n-k}n\binom{n+k}{n-k} 
\]
\[
A_n=\binom{2n}n\sum_k\binom{2k}k\binom{2n-2k}{n-k}\binom{2n-k}k\binom{n+k}{%
n-k} 
\]
\[
A_n=\sum_k\binom nk\binom{n+k}n^2\binom{2k}n\binom{2n}{n-k} 
\]
\[
A_n=\binom{2n}n^2\sum_k\binom nk^2\binom{2n-k}n\binom{n+k}{n-k}\binom{2n}{2k}%
^{-1} 
\]
\[
A_n=\binom{2n}n\sum_k(-1)^k\binom nk\binom{2k}k\binom{2n-2k}{n-k}\binom{n+2k}%
n 
\]
\[
A_n=\binom{2n}n\sum_k(-1)^{n+k+1}\binom nk\binom{2k}k\binom{2n-k}n\binom{%
3n-2k}{n-k} 
\]
\[
A_n=\binom{2n}n\sum_k(-1)^k\binom nk\binom{2k}k\binom{2n-2k}{n-k}\binom{3n-2k%
}n 
\]
\[
A_n=\binom{2n}n\sum_k(-1)^{n+k}\binom{n+k}n^3\binom{3n+1}{n-k} 
\]
\[
A_n=\binom{2n}n\sum_k\binom{n+k}n^3\binom{3n+1}{n-2k} 
\]
\[
A_n=\frac{(3n+1)!}{n!^3(2n+1)}\sum_k(-1)^{n+k}\binom nk\binom{n+k}n^3\binom{%
2n+k+1}k^{-1} 
\]
\[
A_n=\binom{2n}n\sum_{i,j}\binom ni^2\binom nj^2\binom{i+j}n 
\]
\[
A_n=\sum_{i,j}\binom ni\binom nj\binom{2n-i}n\binom{n+j}n\binom ij\binom{2n}%
i 
\]
\[
A_n=\sum_{i,j}\binom ni^2\binom nj\binom{n+j}n\binom{n+i-j}n\binom{2n}{n-j} 
\]
\[
A_n=\sum_{i,j}\binom ni^2\binom nj\binom{n+j}n\binom{i+j}n\binom{2n}{n-j} 
\]
\[
A_n=\sum_{i,j}\binom ni^2\binom nj\binom{i+j}j\binom{2n-i}n\binom{n+i}{n-j} 
\]
\[
A_n=\sum_{i,j}\binom ni^2\binom nj\binom{i+j}j\binom{2i}n\binom{n+i}{n-j} 
\]
\[
A_n=\sum_{i,j}\binom ni^2\binom nj\binom{n+i-j}{n-j}\binom{2n+j}{2n}\binom{%
n+i}{n+j} 
\]
\[
A_n=\binom{2n}n\sum_{i,j}\binom ni^2\binom nj^2\binom{n+i-j}n 
\]
\[
A_n=\binom{2n}n\sum_{i,j}\binom ni^2\binom nj\binom ij\binom{n+i-j}{n-j} 
\]
\[
A_n=\binom{2n}n\sum_{i,j}\binom ni^2\binom nj\binom ij\binom{n+j}n 
\]
\[
A_n=\binom{2n}n\sum_{i,j}\binom ni^2\binom nj\binom ij\binom{2n-i}n 
\]
\[
A_n=\binom{2n}n\sum_{i,j}\binom ni^2\binom nj\binom{2n-i}n\binom i{n-j} 
\]
\[
A_n=\binom{2n}n\sum_{i,j}\binom ni^2\binom nj\binom{2i}n\binom i{n-j} 
\]
\[
A_n=\binom{2n}n\sum_{i.j}\binom ni^2\binom nj\binom{i+j}j\binom i{n-j} 
\]
\[
A_n=\binom{2n}n\sum_{i.j}\binom ni^2\binom{2n-i}n\binom{n+i-j}{n-j}\binom{n+i%
}j 
\]
\[
A_n=\binom{2n}n\sum_{i.j}\binom ni\binom nj\binom{n+j}n\binom ij\binom
n{i-j} 
\]
\[
A_n=\binom{2n}n\sum_{i.j}\binom ni\binom nj\binom{2n-j}n\binom ij\binom{2j}i 
\]
\[
A_n=\binom{2n}n\sum_{i.j}\binom ni\binom ij\binom{2n-j}n\binom{2j}j\binom{2j%
}{i-j} 
\]
\[
A_n=\binom{2n}n\sum_{i,j}\binom ni^2\binom nj\binom ij\binom{2i}n 
\]
\[
A_n=\binom{2n}n\sum_{i,n}(-1)^{i+j}\binom{2n-i}n^2\binom{2n-j}n\binom{n+j}n%
\binom{2n+1}{i-j} 
\]
\[
A_n=\binom{2n}n\sum_{i,j}(-1)^{i+j}\binom{2n-i}n^2\binom{n+j}n^2\binom{3n+1}{%
i-j} 
\]
\[
A_n=\binom{2n}n\sum_{i,j}(-1)^j\binom ni^2\binom{i+j}j\binom{n+i}n\binom{3n+1%
}{n-j} 
\]
\[
A_n=\sum_{i,j}(-1)^j\binom ni^2\binom{2n-i}n\binom{n+i}n\binom{n+j}n\binom{%
3n+1}{n-j} 
\]
\[
A_n=\sum_{i,j}(-1)^{n+j}\binom ni^2\binom{i+j}n\binom{n+j}n\binom{n+i+j}n%
\binom{3n+1}{n-j} 
\]
\[
A_n=\sum_{i,j}(-1)^{n+j}\binom ni^2\binom{n+j}n^2\binom{n+j}n\binom{3n+1}{%
n-2j} 
\]
\[
A_n=\binom{2n}n\sum_{i,j}(-1)^{i+j}\binom ni\binom nj\binom{n+j}n\binom{2n-j}%
n\binom n{2j-i} 
\]
\[
A_n=\binom{2n}n\sum_{i,j}(-1)^{i+j}\binom ni\binom nj\binom{2n-i}n\binom{2n-j%
}n\binom{2n-i-j}n 
\]
\[
"A"_n=\binom{2n}n\sum_k(n-2k)\binom nk^6 
\]
\[
\]

\textbf{27.} 
\[
A_n=\sum_{i,j}\binom ni^2\binom nj^2\binom{i+j}n\binom{2n-i}n 
\]
\[
A_n=\sum_{i,j}\binom ni^2\binom nj^2\binom{n+j}n\binom{n+i-j}n 
\]
\[
A_n=\sum_{i,j}\binom ni^2\binom nj^2\binom{2n-i-j}{n-j}\binom{2i+j}{n+j} 
\]
\[
A_n=\sum_{i,j}\binom ni^2\binom nj\binom ij\binom{n+j}n\binom{2n-i}n 
\]
\[
A_n=\sum_{i.j}\binom ni^2\binom nj\binom ij\binom{n+i-j}{n-j}\binom{2n-i}n 
\]
\[
A_n=\sum_{i,j}\binom ni^2\binom nj\binom ij\binom{n+j}n\binom{n+i-j}{n-j} 
\]
\[
A_n=\sum_{i,j}\binom ni^2\binom nj\binom{i+j}j\binom{2n-i-j}n\binom{n+i}{n-j}
\]
\[
A_n=\sum_{i,j}\binom ni^2\binom nj\binom{2n-i}n\binom{n+j}n\binom ij 
\]
\[
A_n=\sum_{i,j}\binom ni^2\binom nj\binom{2n-i}n\binom{i+j}j\binom i{n-j} 
\]
\[
A_n=\sum_{i.j}\binom ni\binom nj\binom{2n-i}n\binom{n+j}n\binom{i+j}n\binom
n{i-j} 
\]
\[
A_n=\sum_{i,j}\binom ni^2\binom nj\binom{n+i}n\binom{n+j}n\binom{n-i}j 
\]
\[
A_n=\sum_{i,j}\binom ni^2\binom nj^2\binom{2n-i-j}n\binom{n+j}n 
\]
\[
\]

\textbf{28. } 
\[
A_n=\sum_{i,j}\binom ni^2\binom nj^2\binom{i+j}n^2 
\]
\[
A_n=\sum_{i,j}\binom ni^2\binom nj^2\binom{i+j}j\binom{2n-i-j}{n-j} 
\]
\[
A_n=\sum_{i,j}\binom ni^2\binom nj\binom ij\binom{i+j}j\binom{2n-j}n 
\]
\[
A_n=\sum_{i,j}\binom ni^2\binom nj\binom ij\binom{n+i}n\binom{2n-i-j}{n-j} 
\]
\[
A_n=\sum_{i,j}\binom ni^2\binom nj\binom{i+j}n\binom{n+j-i}n\binom{2n-i}n 
\]
\[
A_n=\sum_{i,j}\binom ni^2\binom nj\binom ij\binom{i+j}j\binom{n+j}n 
\]
\[
A_n=\sum_{i,j}\binom ni^2\binom nj\binom{n+j}n\binom{2j}j\binom{i+j}{2j} 
\]
\[
A_n=\sum_{i,j}\binom ni^2\binom nj\binom{n+j}n\binom{n+i-j}{n-j}\binom
i{n-j} 
\]
\[
A_n=\sum_{i,j}\binom ni\binom nj\binom{n+i}n\binom{n+j}n\binom{2n-i-j}%
n\binom n{i+j} 
\]
\[
A_n=\sum_{i,j}\binom ni\binom nj\binom{n+i}n\binom{i+j}j\binom{n+i-j}i\binom
n{i+j} 
\]
\[
A_n=\sum_{i,j}\binom ni\binom nj\binom{2n-i}n\binom{n+j}n\binom{n+i-j}%
n\binom n{i-j} 
\]
\[
A_n=\sum_{i,j}\binom ni\binom nj\binom{2n-i}n\binom{n+j}n\binom{2j}j\binom
j{i-j} 
\]
\[
A_n=\sum_{i,j}\binom ni\binom nj\binom{2n-i}n\binom{2n-j}n\binom{2j}j\binom
ij 
\]
\[
A_n=\sum_{i,j}\binom ni\binom nj\binom{2n-i}n\binom{n+j}n\binom{2j}i\binom
ij 
\]
\[
A_n=\sum_{i,j}\binom ni^2\binom nj\binom{n+j}n\binom{n+i-j}n\binom{2i}{n+j} 
\]
\[
A_n=\sum_{i,j}\binom ni^2\binom nj\binom{n+j}n\binom{n+i-j}{n-j}\binom{2n-i}{%
n+j} 
\]
\[
\]

\textbf{29.} 
\[
A_n=\binom{2n}n\sum_k\binom nk\binom{n+k}n\binom{2k}k\binom{n+k}{n-k} 
\]
\[
A_n=\binom{2n}n\sum_k\binom nk^2\binom{n+k}n^2 
\]
\[
A_n=\binom{2n}n^2\sum_k\binom nk^3\binom{2n-k}n\binom{2n}k^{-1} 
\]
\[
A_n=\sum_{i,j}\binom ni^2\binom nj^2\binom{i+j}j\binom{n+i+j}n 
\]
\[
A_n=\sum_{i,j}\binom ni^2\binom nj\binom ij\binom{n+i}n\binom{2n+i-j}{n-j} 
\]
\[
A_n=\sum_{i,j}\binom ni^2\binom nj^2\binom{n+i-j}n\binom{2n+i-j}n 
\]
\[
A_n=\sum_{i,j}\binom ni^2\binom nj\binom ij\binom{i+j}n\binom{n+i-j}{n-j} 
\]
\[
A_n=\sum_{i,j}\binom ni^2\binom nj\binom{i+j}j\binom{n+i}n\binom{n+i}{n-j} 
\]
\[
A_n=\binom{2n}n\sum_{i,j}\binom ni^2\binom nj\binom{i+j}n\binom{n+i}{n-j} 
\]
\[
A_n=\binom{2n}n\sum_{i,j}\binom ni^2\binom nj\binom{2n-i-j}{n-j}\binom{n+i}{%
n-j} 
\]
\[
A_n=\binom{2n}n\sum_{i,j}\binom ni^2\binom nj\binom{n+j}n\binom{n+i}{n+j} 
\]
\[
A_n=\binom{2n}n\sum_{i,j}\binom ni^2\binom nj\binom{i+j}j\binom{n+i}{n+j} 
\]
\[
A_n=\sum_{i,j}\binom ni\binom nj\binom{2n-i}n\binom{2n-j}n\binom ij\binom{2n}%
i 
\]
\[
A_n=\sum_{i,j}\binom ni^2\binom nj\binom{n+j}n\binom{n+i-j}{n-j}\binom{2n}{%
n-j} 
\]
\[
A_n=\sum_{i,j}\binom ni^2\binom nj\binom{n+j}n\binom{3n-i-j}{n-j}\binom{2n-i%
}{n+j} 
\]
\[
A_n=\binom{2n}n\sum_{i,j}\binom ni^2\binom nj^2\binom{i+j}n 
\]
\[
A_n=\binom{2n}n\sum_{i,j}\binom ni^2\binom nj^2\binom{2n+i-j}{i-j} 
\]
\[
A_n=\sum_{i,j}\binom ni^2\binom nj^2\binom ij\binom{n+i}n 
\]
\[
A_n=\binom{2n}n\sum_{i,j}\binom ni\binom{n+i}n\binom ij^3 
\]
\[
A_n=\binom{2n}n\sum_{i,j}\binom ni\binom ij\binom{n+j}n\binom{2j}j\binom{2j}{%
i-j} 
\]
\[
A_n=\binom{2n}n\sum_{i,j}\binom ni^2\binom nj\binom{n+i-j}{n-j}\binom{2i}{n+j%
} 
\]
\[
A_n=\sum_{i,j}\binom ni^2\binom nj^2\binom{i+j}n\binom{n+i+j}n 
\]
\[
A_n=\sum_{i,j}\binom ni^2\binom nj^2\binom{i+j}n\binom{2n-j}n 
\]
\[
A_n=\sum_{i,j}\binom ni\binom nj\binom{2n-i}n\binom{i+j}n^2\binom{2n}{i+j} 
\]
\[
A_n=\binom{2n}n\sum_{i,j}\binom ni\binom nj\binom{n+i}n\binom{2j}j\binom ij 
\]
\[
A_n=\binom{2n}n\sum_{i,j}\binom ni\binom ij^3\binom{n+i}n 
\]
\[
A_n=\binom{2n}n\sum_{i,j}(-1)^{i+j+n}\binom ni\binom nj\binom{n+i}n\binom{n+j%
}n\binom{i+j}j 
\]
\[
A_n=\binom{2n}n\sum_{i,j}(-1)^{i+j}\binom ni\binom nj\binom{n+i}n\binom{n+j}n%
\binom{n+i+j}n 
\]
\[
A_n=\binom{2n}n\sum_{i,j}(-1)^{i+j}\binom ni\binom nj\binom{n+j}n\binom{i+j}%
n^2 
\]

\[
A_n=\binom{2n}n\sum_{i,j}(-1)^{i+j}\binom{n+i}n^2\binom{n+j}n^2\binom{2n+1}{%
i-j} 
\]
\[
A_n=\binom{2n}n\sum_{i,j}\binom ni^2\binom nj^2\binom{2n+i-j}{2n} 
\]
\[
A_n=\binom{2n}n\sum_{i,j}\binom ni^2\binom nj\binom{n+i}n\binom ij 
\]
\[
A_n=\binom{2n}n\sum_{i,j}\binom ni^2\binom nj\binom{n+j}n\binom{2n-i}{n+j} 
\]
\[
A_n=\binom{2n}n\sum_{i,j}\binom ni^2\binom nj\binom{n+i-j}n\binom{n+i}j 
\]
\[
A_n=\binom{2n}n\sum_{i,j}\binom ni\binom nj\binom{n+i}n\binom{n+j}n\binom
n{i+j} 
\]
\[
A_n=\binom{2n}n\sum_{i,j}(-1)^{i+j}\binom{i+j}n\binom{n+i}n^2\binom{n+j}%
n\binom n{i+j} 
\]
\[
A_n=\binom{2n}n\sum_{i,j}(-1)^{i+j}\binom{i+j}n\binom{2n-i}n^2\binom{n+j}%
n\binom n{i+j} 
\]
\[
A_n=\sum_{i,j}(-1)^j\binom ni^2\binom{2n-i}n^2\binom{n+j}n\binom{3n+1}{n-j} 
\]
\[
A_n=\sum_{i,j}\binom ni\binom nj\binom{n+i}n\binom{n+j}n\binom{i+j}j\binom
n{i+j} 
\]
\[
A_n=\sum_{i,j}\binom ni\binom nj\binom{n+i}n\binom{n+j}n\binom{n+i-j}n\binom
n{i+j} 
\]
\[
A_n=\binom{2n}n\sum_{i,j}\binom ni\binom nj\binom{2n-i}n\binom{n+j}n\binom
n{i-j} 
\]
\[
A_n=\binom{2n}n\sum_{i,j}\binom ni\binom nj\binom{2n-j}n\binom ij\binom
n{i-j} 
\]
\[
A_n=\sum_{i,j}\binom ni\binom nj\binom{n+j}n^3\binom n{i-j} 
\]
\[
A_n=\binom{2n}n\sum_{i,j}\binom ni\binom nj\binom{n+j}n\binom{n+i-j}n\binom
n{i-j} 
\]
\[
A_n=\binom{2n}n\sum_{i,j}\binom ni\binom nj\binom{n+j}n\binom ij\binom{2j}i 
\]
\[
A_n=\binom{2n}n\sum_{i,j}\binom ni\binom nj\binom{n+j}n^2\binom n{2j-i} 
\]
\[
A_n=\binom{2n}n\sum_{i,j}\binom ni\binom nj\binom{2n-j}n^2\binom n{2j-i} 
\]
\[
A_n=\binom{2n}n\sum_{i,j}(-1)^{i+j}\binom{2n-i}n^2\binom{n+j}n^2\binom{2n+1}{%
i-j} 
\]
\[
"A_n"=(2n+1)^{-1}\binom{2n}n^{-1}\sum_k\binom{n+k}n^4\binom{2n+k+1}k^{-2} 
\]

\textbf{13*.} 
\[
A_n=432^n\binom{2n}n\sum_k\binom{-1/6}k^2\binom{-5/6}{n-k}^2 
\]
\[
A_n=432^n\binom{2n}n\sum_k(-1)^{n+k}\binom nk\binom{n+k}n\binom{-1/6}k\binom{%
-5/6}k 
\]
\[
\]

\textbf{13**.} 
\[
A_n=\binom{2n}n\sum_k\frac{(6k)!}{(3k)!(2k)!k!}\frac{(6n-6k)!}{%
(3n-3k)!(2n-2k)!(n-k)!} 
\]
\[
A_n=432^n\binom{2n}n\sum_k\binom{-1/6}k\binom{-5/6}k\binom{-1/6}{n-k}\binom{%
-5/6}{n-k} 
\]
\[
A_n=432^n\binom{2n}n\sum_k(-1)^{n+k}\binom nk\binom{2k}n\binom{-1/6}k\binom{%
-5/6}k 
\]
\[
\]

\textbf{32.} 
\[
A_n^{\prime }=\sum_{i,j}\binom ni^2\binom nj^2\binom{n+i}n\binom{n+j}n\binom{%
i+j}n 
\]
\[
A_n^{\prime }=\sum_{i,j}\binom ni^2\binom nj\binom{n+i-j}n\binom{n+i-j}{n-j}%
\binom{2n-j}n\binom{n+i}j 
\]
\[
A_n^{\prime }=\sum_{i,j}\binom ni^2\binom nj\binom{n+j}n\binom{n+i-j}{n-j}%
\binom{2n-j}n\binom{n+i}{n+j} 
\]
\[
A_n^{\prime }=\sum_{i,j}\binom ni\binom nj\binom{n+j}n^2\binom{n+i}n\binom{%
2n-i}n\binom n{i-j} 
\]
\[
A_n^{\prime }=\sum_{i,j}(-1)^{i+j}\binom{2n-i}n^3\binom{n+j}n^3\binom{3n+1}{%
i-j} 
\]
\[
A_n^{\prime }=\sum_{i,j}\binom ni^2\binom nj\binom{n+j}n\binom{2n-j}n\binom{%
2n-i-j}{n-j}\binom{2n-i}{n+j} 
\]
\[
A_n^{\prime }=\sum_{i,j}(-1)^j\binom ni^2\binom{n+i}n\binom{n+j}n^2\binom{%
n+i+j}n\binom{3n+1}{n-j} 
\]
\[
"A_n^{\prime }"=\sum_k(n-2k)\binom nk^4\binom{n+k}n^2\binom{2n-k}n^2 
\]
\[
\]

\textbf{33.} 
\[
A_n=\binom{2n}n^2\sum_k\binom{2n}k^3 
\]
\[
A_n=64^n\binom{2n}n^2\sum_k(-1)^k\binom{2n+k}{2n}^3 
\]
\[
A_n=\binom{2n}n^2\sum_k(-1)^k\binom nk\binom{2k}k\binom{2n-k}n\binom{4n-2k}{%
2n-k}\binom{n+k}n^{-1} 
\]
\[
A_n=\binom{2n}n^3\sum_k(-1)^k\binom nk^2\binom{2k}k\binom{4n-2k}{2n-k}\binom{%
2n}k^{-1}\binom{n+k}n^{-1} 
\]
\[
A_n=\binom{2n}n^2\sum_{i,j}\binom ni\binom nj\binom{2n}{i+j}^2 
\]
\[
"A_n"=64^n\binom{2n}n^2\sum_kk^{-3}\binom{2n}k^{-3} 
\]
\[
"A_n"=16^{-n}\binom{2n}n^3\sum_k(-1)^{n+k}\binom nk\binom{2k}k\binom{2n-k}n%
\binom{4n-2k}{2n-k}\binom{n+k}n^{-1}\binom{2n}{2k}^{-1} 
\]
\[
\]

\textbf{\ }

\textbf{36.} 
\[
A_n=\binom{2n}n^2\sum_k\binom nk\binom{2k}k\binom{2n-2k}{n-k} 
\]
\[
A_n=2^{-n}\binom{2n}n^4\sum_k\binom nk^3\binom{2n}{2k}^{-1} 
\]
\[
A_n=2^{-n}\binom{2n}n\sum_k\binom{n+k}n\binom{2k}k^2\binom{2n-2k}{n-k}^2%
\binom{2n}{n-k}\binom nk^{-2} 
\]
\[
A_n=2^{-n}\binom{2n}n^2\sum_k\binom nk^{-1}\binom{2k}k^2\binom{2n-2k}{n-k}^2 
\]
\[
\]

\textbf{38.} 
\[
A_n=\binom{2n}n\binom{4n}{2n}\sum_k\binom nk\binom{2k}k\binom{2n-2k}{n-k} 
\]
\[
A_n=4^{-n}\binom{2n}n^2\binom{4n}{2n}\sum_k\binom nk^3\binom{2n}{2k}^{-1} 
\]
\[
\]

\textbf{39.} 
\[
A_n^{\prime }=\binom{2n}n^2\sum_k\binom nk^2\binom{2k}k\binom{2n-2k}{n-k} 
\]
\[
A_n^{\prime }=\sum_k\binom{n+k}n\binom{2n-k}n\binom{2k}k\binom{2n-2k}{n-k}%
\binom{2n}k\binom{2n}{n-k} 
\]
\[
A_n^{\prime }=\binom{2n}n\sum_k\binom{2n}{2k}\binom{2k}k^2\binom{2n-2k}{n-k}%
^2 
\]
\[
\]

\textbf{40.} 
\[
A_n^{\prime }=\binom{2n}n^2\sum_k\binom{2k}k^2\binom{2n-2k}{n-k}^2 
\]
\[
A_n^{\prime }=\binom{2n}n\sum_k\binom{n+k}n\binom{2n-k}n\binom{2k}k\binom{%
2n-2k}{n-k}\binom{2n}k\binom{2n}{n-k}\binom{2n}{2k}^{-1} 
\]
\[
\]

\textbf{42.} 
\[
A_n=\binom{2n}n\sum_{i,j}\binom nk^2\binom{2k}n^2 
\]
\[
A_n=\binom{2n}n\sum_k\binom nk\binom{2k}k\binom{2k}n\binom k{n-k} 
\]
\[
A_n=\binom{2n}n\sum_{i,j}\binom ni^2\binom nj\binom ij\binom{n+i-j}n 
\]
\[
A_n=\binom{2n}n\sum_{i,j}\binom ni^2\binom nj\binom ij\binom{i+j}j 
\]
\[
A_n=\binom{2n}n\sum_{i,j}\binom ni^2\binom nj\binom{i+j}n\binom i{n-j} 
\]
\[
A_n=\binom{2n}n\sum_{i,j}\binom ni^2\binom nj\binom{n-i+j}j\binom i{n-j} 
\]
\[
A_n=\binom{2n}n\sum_{i,j}\binom ni\binom nj\binom{2n-i}n\binom ij\binom
n{i-j} 
\]
\[
A_n=\binom{2n}n\sum_{i,j}\binom ni\binom nj\binom{2n-i}n\binom ij\binom{2j}i 
\]
\[
A_n=\binom{2n}n\sum_{i,j}\binom ni^2\binom nj\binom{2j}j\binom{i+j}{2j} 
\]
\[
A_n=\binom{2n}n\sum_{i,j}\binom ni^2\binom nj\binom{n+j}n\binom{2i}{n+j} 
\]
\[
\]

\textbf{44.} 
\[
A_n^{\prime }=\binom{2n}n^2\sum_k\binom nk^2\binom{n+k}n^2 
\]
\[
A_n^{\prime }=\sum_{i,j}(-1)^{i+j}\binom{2n-i}n^3\binom{n+j}n^2\binom{2n+j}n%
\binom{3n+1}{i-j} 
\]
\[
A_n^{\prime }=\sum_{i,j}(-1)^{i+j}\binom{2n-i}n^2\binom{n+j}n^3\binom{3n-i}n%
\binom{3n+1}{i-j} 
\]

\[
A_n^{\prime }=\binom{2n}n\sum_{i,j}\binom ni^2\binom nj^2\binom{i+j}n^2 
\]
\[
"A_n^{\prime }"=\binom{2n}n^2\sum_kk^{-4}\binom nk^{-2}\binom{n+k}n^{-2} 
\]
\[
"A_n^{\prime }"=(2n+1)^{-2}\sum_k\binom{n+k}n^4\binom{2n+k+1}k^{-2} 
\]
\[
\]

\textbf{45.} 
\[
A_n=\binom 2n^2\sum_k\binom nk^3 
\]
\[
A_n=\sum_k\binom nk\binom{n+k}n\binom{2n-k}n\binom{2n}k\binom{2n}{n-k} 
\]
\[
A_n=\binom{2n}n\sum_k\binom nk^2\binom{2n-k}n\binom{2n}k 
\]
\[
A_n=2^{-n}\binom{2n}n^2\sum_k\binom{2k}k\binom{2n-2k}{n-k}\binom{2k}n 
\]
\[
A_n=\binom{2n}n\sum_k\binom nk\binom{n+k}n\binom{2k}n\binom{2n}{n-k} 
\]
\[
A_n=\binom{2n}n\sum_k\binom nk\binom{2n-k}n\binom{2k}n\binom{2n}k 
\]
\[
A_n=\binom{2n}n^2\sum_k\binom nk\binom{2k}k\binom k{n-k} 
\]
\[
A_n=8^n\binom{2n}n^4\sum_k(-1)^k\binom nk^2\binom{2n-k}n\binom{2n}k^{-2} 
\]
\[
A_n=8^n\binom{2n}n^5\sum_k(-1)^k\binom nk^3\binom{2n}k^{-3} 
\]
\[
A_n=8^n\binom{2n}n^2\sum_k(-1)^{n-k}\binom{n+k}n^3 
\]
\[
A_n=\binom{2n}n\sum_{i,j}\binom ni^2\binom nj\binom ij\binom{2n}j 
\]
\[
A_n=\sum_{i,j}\binom ni\binom nj\binom{2n-i}n\binom ij\binom{2n}i\binom{2n}j 
\]
\[
A_n=\binom{2n}n\sum_{i,j}\binom ni\binom nj\binom{n+j}n\binom{2j}n\binom
n{i-j} 
\]
\[
A_n=\binom{2n}n\sum_{i,j}\binom ni^2\binom nj\binom{n-j}i\binom{n+i}{n-j} 
\]
\[
A_n=\binom{2n}n^2\sum_{i,j}\binom ni\binom nj\binom ij\binom j{i-j} 
\]
\[
A_n=\binom{2n}n^2\sum_{i,j}(-1)^{i+j}\binom ni\binom nj^2\binom n{2j-i} 
\]
\[
A_n=\sum_{i,j}\binom ni^2\binom nj\binom{i+j}j\binom{2n-i-j}{n-j}\binom{2n}{%
i+j} 
\]
\[
A_n=\binom{2n}n\sum_{i,j}\binom ni^2\binom nj\binom ij\binom{2n}{n+i} 
\]
\[
A_n=\binom{2n}n\sum_{i,j}(-1)^{i+j}\binom ni\binom nj\binom ij\binom{3n-i}n%
\binom{2n-j}n 
\]
\[
A_n=\binom{2n}n^2\sum_{i,j}\binom ni^2\binom nj\binom{n-2i}j 
\]
\[
A_n=\binom{2n}n\sum_{i,j}(-1)^{i+j}\binom{2n-i}n^2\binom{n+j}n\binom{2n}{n-j}%
\binom{2n+1}{i-j} 
\]
\[
A_n=\sum_{i,j}(-1)^{i+j}\binom{2n-i}n\binom{n+j}n^2\binom{2n+j}n\binom{2i}n%
\binom{3n+1}{i-j} 
\]
\[
A_n=\sum_{i,j}(-1)^{i+j}\binom ni\binom{2n-i}n\binom{n+j}n^2\binom{2n+j}n%
\binom{3n+1}{i-j} 
\]
\[
"A_n"=\binom{2n}n^2\sum_k(-1)^kk^{-3}\binom{n+k}n^{-3} 
\]
\[
"A_n"=8^n\binom{2n}n^2\sum_kk^{-3}\binom nk^{-3} 
\]
\[
\]

\textbf{46.} 
\[
A_n=\binom{2n}n\sum_k\frac{(3k)!}{k!^3}\frac{(3n-3k)!}{(n-k)!^3} 
\]
\[
A_n=27^n\binom{2n}n\sum_k(-1)^{n+k}\binom nk\binom{2k}n\binom{-1/3}k\binom{%
-2/3}k 
\]
\[
A_n=27^n\binom{2n}n\sum_{i,j}\binom ni\binom nj\binom{n+i+j}n\binom{-1/3}i%
\binom{-2/3}j 
\]
\[
\]

\textbf{47} 
\[
A_n=\binom{2n}n\sum_k\frac{(6k)!}{(3k)!(2k)!k1}\frac{(6n-6k)!}{%
(3n-3k)!(2n-2k)!(n-k)!} 
\]
\[
A_n=432^n\binom{2n}n\sum_k(-1)^k\binom nk\binom{2k}n\binom{-1/6}k\binom{-5/6}%
k 
\]
\[
\]

\textbf{49.} 
\[
A_n^{\prime }=\frac{(3n)!}{n!^3}\sum_k\binom{2k}k^2\binom{2n-2k}{n-k}^2 
\]
\[
A_n^{\prime }=\binom{2n}n^2\sum_k\binom nk^2\binom{n+k}n\binom{2k}k\binom{%
2n-2k}{n-k}\binom{3n}{n+k}\binom{2n}k^{-1}\binom{2n}{2k}^{-1} 
\]
\[
\]

\textbf{50.} 
\[
A_n^{\prime }=\frac{(3n)!}{n!^3}\sum_{i+j+k+m=n}\left\{ \frac{n!}{i!j!k!m!}%
\right\} ^2 
\]
\[
A_n^{\prime }=\frac{(3n)!}{n!^3}\sum_k\binom nk^2\binom{2k}k\binom{2n-2k}{n-k%
} 
\]
\[
A_n^{\prime }=4^{-n}\frac{(3n)!}{n!^3}\binom{2n}n^3\sum_k\binom nk^4\binom{2n%
}{2k}^{-3} 
\]
\[
A_n^{\prime }=\binom{2n}n\sum_k(-1)^{n+k}\binom nk\binom{n+k}n\binom{2k}k%
\binom{2n-2k}{n-k}\binom{2k}n\binom{3n}{n+k}\binom{2n}k^{-1} 
\]
\[
\]

\textbf{53.} 
\[
A_n=\frac{(3n)!}{n!^3}\sum_k\binom nk^2\binom{n+k}n^2 
\]
\[
A_n=\binom{2n}n^2\sum_k\binom nk^3\binom{n+k}n\binom{2n-k}n\binom{3n}{n+k}%
\binom{2n}k^{-2} 
\]
\[
\]

\textbf{55.} 
\[
A_n=\binom{2n}n^2\sum_k\binom nk^2\binom{2n}{2k} 
\]
\[
A_n=\binom{2n}n\sum_k\binom{2k}k\binom{2n-2k}{n-k}\binom{2n}{2k}^2 
\]
\[
A_n=\binom{2n}n^2\sum_k(-1)^{n+k}\binom nk\binom{n+k}n\binom{2n}k\binom{2n}{%
n-k}\binom{2n}{2k}^{-1} 
\]
\[
A_n=\binom{2n}n^2\sum_k\binom nk\binom{2n-k}n\binom{2n}k^2\binom{2n}{2k}%
^{-1} 
\]
\[
A_n=\binom{2n}n^2\sum_k(-1)^{n+k}\binom{n+k}k\binom{2n-2k}{n-k}\binom{2n}{n-k%
}\binom{2n+2k}{n+k}\binom{2n}{2k}^{-1} 
\]
\[
A_n=\binom{2n}n\sum_k(-1)^{n+k}\binom{n+k}k\binom{2n-2k}{n-k}\binom{2n}{n-k}%
\binom{2n+2k}{n+k} 
\]
\[
A_n=\binom{2n}n^2\sum_k(-1)^k\binom nk\binom{2k}k\binom{4n-2k}{2n-k} 
\]
\[
A_n=\binom{2n}n^2\sum_k(-1)^k\binom{2k}k\binom{2n-2k}{n-k}\binom{2n}k 
\]
\[
A_n=\binom{2n}n\sum_k(-1)^k\binom{2k}k\binom{2n}k^2\binom{4n-2k}{2n-k} 
\]
\[
A_n=\binom{2n}n^3\sum_k(-1)^k\binom nk^2\binom{2n}k\binom{2n}{2k}^{-1} 
\]
\[
A_n=\binom{2n}n^2\sum_k(-1)^k\binom nk\binom{2n-k}n\binom{2n}k^2\binom{2n}{2k%
}^{-1} 
\]
\[
A_n=\binom{2n}n\sum_k(-1)^k\binom{n+k}k^2\binom{2n}{n+k}^3\binom{2n}{2k}%
^{-1} 
\]
\[
A_n=4^{-n}\binom{2n}n^2\sum_k(-1)^k\binom nk\binom{2k}k\binom{2n-k}n\binom{%
4n-2k}{2n-k}\binom{2n}{2k}^{-1} 
\]
\[
A_n=\binom{2n}n^3\sum_k(-1)^k\binom nk\binom{2k}k\binom{2n-2k}{n-k}\binom{n+k%
}n^{-1} 
\]
\[
A_n=\binom{2n}n\sum_{i,j}\binom ni\binom nj\binom{2i}n\binom{2j}n\binom{2n}{%
i+j} 
\]
\[
A_n=\binom{2n}n^2\sum_{i,j}\binom ni\binom nj^2\binom n{2j-i} 
\]
\[
A_n=\binom{2n}n^2\sum_{i,j}(-1)^{i+j}\binom ni\binom nj\binom{2n}{2j}\binom
n{2j-i} 
\]
\[
"A_n"=16^{-n}\binom{2n}n^3\sum_k(-1)^{n+k}\binom nk\binom{2k}k\binom{4n-2k}{%
2n-k}\binom{2n}{2k}^{-1} 
\]
\[
\]

\textbf{58.} 
\[
A_n=\binom{2n}n^2\sum_k\binom nk^2\binom{2k}k 
\]
\[
A_n=\binom{2n}n^2\sum_k\binom nk^2\binom{3k}{2n} 
\]
\[
A_n=3^{1-n}\binom{2n}n^3\sum_{k=0}^{[n/3]}(-1)^k\frac{n-2k}{2n-3k}\binom nk^2%
\binom{3n-3k}{2n}\binom{2n}{3k}^{-1} 
\]
\[
A_n=\binom{2n}n\sum_{i,j}\binom ni\binom nj\binom ij\binom{n+j}i\binom{3j}{2n%
} 
\]
\[
A_n=\binom{2n}n^2\sum_{i,j}\binom ni\binom ij^3 
\]
\[
A_n=\binom{2n}n^2\sum_{i,j}(-1)^{i+j}\binom ni\binom nj\binom{2j}j\binom
n{2j-i} 
\]
\[
A_n=\binom{2n}n^2\sum_{i,j}(-1)^{i+j}\binom ni\binom nj\binom{3j}{2n}\binom
n{2j-i} 
\]
\[
A_n=\binom{2n}n^2\sum_{i,j}(-1)^{i+j}\binom ni\binom nj\binom{i+j}n^2 
\]
\[
"A_n"=9^n\binom{2n}n^2\sum_k\binom{2k}k^2\binom{n+k}n\binom{n+k}{n-k}\binom
nk^{-1} 
\]
\[
\]

\textbf{59.} 
\[
A_n=\sum_k\binom nk\binom{2n-k}n\binom{n+k}n\binom{2k}k\binom{2n-2k}{n-k} 
\]
\[
"A"_n=\sum_k(n-2k)\binom nk^3\binom{n+k}n\binom{2n-k}n\binom{2k}n\binom{2n-2k%
}n 
\]
\[
\]

\textbf{73.} 
\[
A_n=\binom{2n}n\sum_k\frac{(3k)!}{k!^3}\frac{(3n-3k)!}{(n-k)!^3}\binom{2n}%
k\binom nk^{-1} 
\]
\[
"A"_n=\binom{2n}n^2\sum_k(n-2k)\binom nk^3\binom{3k}n\binom{3n-3k}n 
\]
\[
\]

\textbf{83.} 
\[
A_n=\binom{2n}n\sum_k\binom nk\binom{2n}{2k}^{-1}\frac{(4k)!}{(2k)!k!^2}%
\frac{(4n-4k)!}{(2n-2k)!(n-k)!^2} 
\]
\[
A_n=2\binom{2n}n^2\sum_{k=[(n+1)/2]}^n\binom nk^3\binom{4k}{2n}\binom{2n}{%
4n-4k}^{-1} 
\]
\[
\]

\textbf{99.} 
\[
A_n=\binom{2n}n^2\sum_k\binom nk^2\binom{2n+k}n 
\]
\[
A_n=\binom{2n}n\sum_k\binom{n+k}k\binom{2n-k}n\binom{2n}k^2 
\]
\[
A_n=\binom{2n}n\sum_k(-1)^{n+k}\binom nk\binom{n+k}n^2\binom{2n+k}n 
\]
\[
A_n=\binom{2n}n^3\sum_k\binom nk^2\binom{2n-k}n\binom{n+k}n^{-1} 
\]
\[
A_n=\binom{2n}n^2\sum_k\binom nk\binom{n+k}n\binom{2n}k 
\]
\[
A_n=\binom{2n}n^2\binom{3n}n\sum_k\binom nk\binom{2n}k^2\binom{3n}{n+k}^{-1} 
\]
\[
A_n=\binom{2n}n^2\sum_k(-1)^k\binom nk\binom{2n+k}n^2 
\]
\[
A_n=\binom{2n}n^2\sum_k(-1)^{n+k}\binom nk\binom{n+k}n\binom{2n+k}{2n} 
\]
\[
A_n=\binom{2n}n\sum_{i,j}\binom ni^2\binom{2n-j}n\binom{n+i-j}{n-j}\binom{%
2n+i}{n+j} 
\]
\[
A_n=\sum_{i,j}\binom ni^2\binom nj\binom{i+j}j\binom{n+i}{n-j} 
\]
\[
A_n=\binom{2n}n^2\sum_{i,j}\binom ni\binom nj\binom{2n-i}n\binom n{i+j} 
\]
\[
A_n=\binom{2n}n^2\sum_{i,j}\binom ni\binom nj\binom{n+i}n\binom n{i-j} 
\]
\[
A_n=\binom{2n}n\sum_{i,j}\binom ni\binom nj\binom{n+j}n\binom{2n+j}n\binom
n{i-j} 
\]
\[
A_n=\sum_{i,j}\binom ni\binom nj\binom{n+i}n\binom{n+j}n\binom{2n+i}n\binom
n{i-j} 
\]
\[
A_n=\sum_{i,j}\binom ni\binom nj\binom{n+i+j}n\binom{i+j}n^2\binom{2n}{i+j} 
\]
\[
A_n=\binom{2n}n^2\sum_{i,j}\binom ni\binom nj\binom{n-i}j\binom{2n}{i+j} 
\]
\[
A_n=\binom{2n}n^2\sum_{i,j}(-1)^{i+j}\binom ni\binom nj\binom{2n+j}n\binom
n{2j-i} 
\]
\[
A_n=\binom{2n}n^2\sum_{i,j}(-1)^{i+j}\binom ni\binom nj\binom{3n-j}n\binom
n{2j-i} 
\]
\[
A_n=\sum_{i,j}(-1)^{i+j}\binom ni\binom nj\binom{i+j}n^3\binom{n+i+j}n 
\]
\[
A_n=\binom{2n}n\sum_{i,j}(-1)^{i+j}\binom ni\binom nj\binom{i+j}n\binom{n+i+j%
}n^2 
\]
\[
A_n=\binom{2n}n\sum_{i,j}\binom ni^2\binom nj\binom{n+i-j}{n-j}\binom{n+2i}{%
n+j} 
\]
\[
A_n=\sum_{i,j}(-1)^{i+j}\binom{2n-i}n\binom{n+j}n^2\binom{2n+j}n\binom{3n-i}n%
\binom{3n+1}{i-j} 
\]
\[
A_n=\sum_{i,j}(-1)^{i+j}\binom{2n-i}n^2\binom{n+j}n\binom{2n+j}n\binom{3n-i}n%
\binom{3n+1}{i-j} 
\]
\[
A_n=\sum_{i,j}(-1)^{i+j}\binom{2n-i}n^3\binom{n+j}n\binom{3n-i}n\binom{3n+1}{%
i-j} 
\]
\[
A_n=\sum_{i,j}(-1)^{i+j}\binom{2n-i}n\binom{n+i}n^2\binom{n+j}n\binom{2n+i}n%
\binom{3n+1}{i-j} 
\]
\[
A_n=\binom{2n}n\sum_{i,j}(-1)^{i+j}\binom{2n+i}n^2\binom{n+j}n\binom{2n+j}n%
\binom{3n+1}{i-j} 
\]
\[
A_n=\binom{2n}n\sum_{i,j}(-1)^{i+j}\binom{2n-i}n\binom{n+j}n\binom{2n+j}n^2%
\binom{3n+1}{i-j} 
\]
\[
A_n=\binom{2n}n\sum_{i,j}(-1)^{i+j}\binom{2n-i}n\binom{n+j}n\binom{3n-i}n^2%
\binom{3n+1}{i-j} 
\]
\[
A_n=\sum_{i,j}(-1)^{i+j}\binom{2n-i}n^2\binom{n+j}n\binom{2n+j}n\binom{3n+1}{%
i-j} 
\]
\[
A_n=\binom{2n}n\sum_{i,j}\binom ni\binom nj\binom{i+j}j\binom{3n-i-j}n\binom{%
2n}{i+j} 
\]
\[
A_n=\sum_{i,j}\binom ni^2\binom nj\binom{2n-j}n\binom{n+i-j}{n-j}\binom{2n+i%
}{n+j} 
\]
\[
A_n=\sum_{i,j}\binom ni^2\binom nj\binom{i+j}j\binom{n+j}n\binom{2n+i}{n+j} 
\]
\[
A_n=\sum_{i,j}(-1)^{n+i+j}\binom{2n-i}n\binom{n+j}n\binom{2n}{i+j}^3\binom{%
2n+1}{i-j} 
\]
\[
A_n=\binom{2n}n\sum_{i,j}(-1)^{i+j}\binom{2n-i}n\binom{n+j}n^2\binom{2n+j}n%
\binom{2n+1}{i-j} 
\]
\[
"A_n"=(-1)^n\binom{2n}n^2\sum_kk^{-3}\binom nk^{-1}\binom{n+k}n^{-1}\binom{2n%
}k^{-1} 
\]
\[
"A_n"=\binom{2n}n^3\sum_k(-1)^{n+k}k^{-3}\binom nk^{-1}\binom{n+k}n^{-1}%
\binom{2n+k}n^{-1} 
\]
\[
"A"_n=\sum_k(n-2k)\binom nk^3\binom{n+k}n\binom{2n-k}n\binom{2n+k}n\binom{%
3n-k}n 
\]
\[
\]

\textbf{101.} 
\[
A_n=\left\{ \sum_k\binom nk^2\binom{n+k}n\right\} ^2 
\]
\[
A_n=\sum_{i,j}\binom ni^2\binom nj\binom{i+j}j\binom{n+j}n\binom{n+i}{n-j} 
\]
\[
A_n=\sum_{i,j}\binom ni^2\binom nj\binom{i+j}j\binom{2n-j}n\binom{n+i}{n-j} 
\]
\[
A_n=\sum_{i,j}\binom ni^2\binom nj^2\binom{n+i-j}i\binom{n+i}j 
\]
\[
A_n=\sum_{i,j}\binom ni^2\binom nj^2\binom{n+i-j}{n-j}\binom{n+i+j}i 
\]
\[
A_n=\sum_{i,j}\binom ni^2\binom nj^2\binom{i+j}j\binom{2n+i-j}i 
\]
\[
A_n=\sum_{i,j}\binom ni^2\binom nj^2\binom{n+i-j}{n-j}\binom{2n-i-j}{n-j} 
\]
\[
A_n=\sum_{i,j}\binom ni^2\binom nj\binom{2n-j}n\binom{2n+i-j}{n-j}\binom{n+i%
}{n+j} 
\]
\[
A_n=\sum_{i,j}\binom ni^2\binom nj\binom{2n-j}n\binom{n+i-j}{n-j}\binom{2n-i%
}{n+j} 
\]
\[
A_n=\sum_{i,j}\binom ni\binom nj\binom{n+j}n^2\binom{i+j}n\binom n{i-j} 
\]
\[
A_n=\sum_{i,j}\binom ni\binom nj\binom{n+i}n\binom{n+j}n\binom{n+i-j}n\binom
n{i-j} 
\]
\[
A_n=\sum_{i,j}\binom ni\binom nj\binom{n+i}n\binom{n+j}n\binom{n+i+j}n\binom
n{i+j} 
\]
\[
A_n=\sum_{i,j}\binom ni\binom nj\binom{n+i}n\binom{2n-j}n\binom{2n-i-j}%
n\binom n{i+j} 
\]
\[
A_n=\sum_{i,j}\binom ni\binom nj\binom{n+i}i\binom{n+j}n\binom ij\binom{2j}i 
\]
\[
A_n=\sum_{i,j}\binom ni\binom nj\binom{n+i}i\binom{n+j}n\binom{2j}j\binom
i{i-j} 
\]
\[
A_n=\sum_{i,j}\binom ni^2\binom nj\binom{i+j}j\binom{n+i}n\binom{2i}{n-j} 
\]
\[
A_n=\sum_{i,j}\binom ni^2\binom nj\binom{n+j}n\binom{2j}j\binom{i+j}{n-j} 
\]
\[
A_n=\binom{2n}n\sum_{i,j}\binom ni^2\binom nj\binom{i+j}n\binom{2i}{n-j} 
\]
\[
A_n=\binom{2n}n\sum_{i,j}\binom ni^2\binom nj\binom ij\binom{2i+j}j 
\]
\[
A_n=\sum_{i,j}\binom ni^2\binom nj\binom{2n-j}n\binom{n+i-j}{n-j}\binom{2i}{%
n+j} 
\]
\[
A_n=\sum_{i,j}\binom ni^2\binom nj\binom{2n-j}n\binom{n+i-j}{n-j}\binom{n+i}%
j 
\]
\textbf{102.} 
\[
A_n=\sum_i\binom ni^3\sum_j\binom nj^2\binom{n+j}n 
\]
\[
A_n=\sum_{i,j}\binom ni^2\binom nj^2\binom{i+j}j\binom{n+i}{n-j} 
\]
\[
A_n=\sum_{i,j}\binom ni^2\binom nj^2\binom{n+i-j}n\binom{n+i}{n-j} 
\]
\[
A_n=\sum_{i,j}\binom ni^2\binom nj\binom{i+j}j\binom{2j}n\binom{n+i}{n-j} 
\]
\[
A_n=\sum_{i,j}\binom ni^2\binom nj^2\binom{2n-i}n\binom{2j}n 
\]
\[
\]

\textbf{109.} 
\[
A_n=\binom{2n}n\sum_k\binom{n+k}n\binom{2n-k}n\binom{2n}k^3\binom nk^{-1} 
\]
\[
A_n=\binom{2n}n^4\sum_k\binom nk^2\binom{2n-k}n\binom{n+k}n^{-2} 
\]
\[
A_n=\binom{2n}n^2\sum_k\binom nk\binom{2n+k}{2n}\binom{3n}{n+k} 
\]
\[
A_n=\binom{2n}n^2\sum_k\binom{2n}k^2\binom{3n-k}n 
\]
\[
A_n=\binom{2n}n^2\sum_k(-1)^{n+k}\binom nk\binom{2n+k}n\binom{3n+k}{2n} 
\]
\[
A_n=\binom{3n}n\sum_k\binom nk\binom{n+k}n\binom{2n+k}n\binom{2n+k}{2n} 
\]
\[
A_n=\frac{(3n)!}{n!^3}\sum_k(-1)^k\binom nk\binom{2n+k}{2n}^2 
\]
\[
A_n=\binom{2n}n^2\sum_k(-1)^k\binom{2n+k}{2n}^2\binom{3n}{n+k} 
\]
\[
A_n=\binom{2n}n\sum_k(-1)^{n+k}\binom nk\binom{2n+k}n^2\binom{3n+k}n 
\]
\[
A_n=\sum_k(-1)^k\binom{n+k}n\binom{2n-k}n\binom{2n+k}{2n}^2\binom{3n}{n+k} 
\]
\[
A_n=\binom{3n}n\sum_k\binom{n+k}n\binom{2n+k}n^2\binom{3n+1}{n-2k} 
\]
\[
A_n=\binom{2n}n\binom{3n}n\sum_{i.j}\binom ni\binom nj\binom n{i+j}\binom{%
2n+i+j}n 
\]
\[
A_n=\sum_{i,j}\binom ni\binom nj\binom{i+j}n^2\binom{n+i+j}n\binom{3n}{i+j} 
\]
\[
"A_n"=\binom{2n}n^2\sum_k\binom nk^2\binom{2n+k}n\binom{3n+k}n\binom{n+k}%
n^{-1} 
\]
\[
\]

\textbf{116.} 
\[
A_n=\binom{2n}n\sum_k4^{n-k}\binom{n+k}n\binom{2k}k^2\binom{2n-2k}{n-k} 
\]
\[
A_n=\binom{2n}n^2\sum_k(-1)^{n+k}\binom nk\binom{2k}k\binom{2n-2k}{n-k}%
\binom{4k}{2n}\binom{2k}n^{-1} 
\]
\[
\]

\textbf{118.} 
\[
A_n=32^{-n}\sum_k\binom nk^{-5}\binom{2k}k^5\binom{2n-2k}{n-k}^5 
\]
\[
"A_n"=32^n\sum_kk^{-5}\binom nk^{-5} 
\]
\[
\]

\textbf{124.} 
\[
A_n=\sum_{i,j}\binom ni^2\binom nj\binom ij\binom{i+j}j\binom{2n-i-j}{n-j} 
\]
\[
A_n=\sum_{i,j}\binom ni^2\binom nj\binom{2j}j\binom{2n-i-j}{n-j}\binom{i+j}{%
2j} 
\]
\[
\]

\textbf{185.} 
\[
A_n=\binom{2n}n\sum_{i,j}\binom ni^2\binom nj\binom ij\binom{i+j}n 
\]
\[
A_n=\binom{2n}n\sum_{i,j}\binom ni\binom nj\binom ij^2\binom n{i-j} 
\]
\[
A_n=\binom{2n}n\sum_{i,j}\binom ni^2\binom nj\binom ij\binom{i+2j}{n+i} 
\]
\[
A_n=\binom{2n}n\sum_{i,j}\binom ni^2\binom nj\binom{2n-i-j}{n-j}\binom{2i}{%
n+i} 
\]
\[
A_n=\binom{2n}n\sum_{i,j}\binom ni^2\binom nj\binom ij\binom{n+j-i}j 
\]
\[
A_n=\binom{2n}n\sum_{i,j}\binom ni\binom nj\binom{i+j}j^2\binom n{i+j} 
\]
\[
A_n=\binom{2n}n\sum_{i,j}\binom ni\binom nj\binom{n+i-j}n\binom{2j}j\binom
j{i-j} 
\]
\[
\]

\textbf{189.} 
\[
A_n=\binom{2n}n\sum_{i,j}\binom ni^2\binom nj^2\binom{i+j}n^2 
\]
\[
\]

\textbf{193.} 
\[
A_n=\sum_{i,j}\binom ni^2\binom nj^2\binom{i+j}j\binom{n+i+j}n 
\]
\[
A_n=\sum_{i,j}\binom ni^2\binom nj\binom{n+i}n\binom{n+i}j\binom{2n-j}{n-j} 
\]
\[
A_n=\sum_{i,j}\binom ni^2\binom nj\binom{2n-j}n\binom{n+i-j}{n-j}\binom{n+i}{%
n+j} 
\]
\[
A_n=\sum_{i,j}\binom ni^2\binom nj\binom{n+i}n\binom{n+j}n\binom{n+i}{n+j} 
\]
\[
A_n=\sum_{i,j}\binom ni\binom nj\binom{n+i}n^2\binom{2n-j}n\binom n{i+j} 
\]
\[
A_n=\sum_{i,j}\binom ni^2\binom nj^2\binom{n+j}n\binom{n+i+j}i 
\]
\[
A_n=\sum_{i,j}\binom ni\binom nj\binom{n+i}n\binom{n+j}n^2\binom n{i-j} 
\]
\[
A_n=\sum_{i,j}\binom ni^2\binom nj^2\binom{2n-j}n\binom{2n+i-j}i 
\]
\[
A_n=\sum_{i,j}\binom ni^2\binom nj\binom{n+i}n\binom{n+j}n\binom{n+i}{n-j} 
\]
\[
A_n=\sum_{i,j}\binom ni^2\binom nj\binom{n+i-j}{n-j}\binom{2n-j}n\binom{n+i}{%
n+j} 
\]
\[
A_n=\sum_{i,j}(-1)^{i+j}\binom ni\binom nj\binom{n+i}n^2\binom{n+j}n\binom{%
i+j}j 
\]
\[
\]

\textbf{194.} 
\[
A_n=\sum_{i,j}\binom ni^2\binom nj^2\binom{i+j}j^2 
\]
\[
A_n=\sum_{i,j}\binom ni^2\binom nj^2\binom{n+i-j}{n-j}\binom{n+i-j}i 
\]
\[
A_n=\sum_{i,j}\binom ni^2\binom nj\binom{i+j}j\binom{i+j}n\binom{2i}{n-j} 
\]
\[
A_n=\sum_{i,j}\binom ni^2\binom nj\binom{i+j}n\binom{2j}j\binom{i+j}{i-j} 
\]
\[
A_n=\sum_{i,j}\binom ni^2\binom nj\binom{i+j}n\binom{2j}j\binom{i+j}{n-j} 
\]
\[
\]

\textbf{195.} 
\[
A_n=\sum_{i,j}\binom ni^2\binom nj^2\binom{i+j}j\binom{n+i}n 
\]
\[
A_n=\sum_{i,j}\binom ni^2\binom nj\binom ij\binom{i+j}j\binom{n+2j}n 
\]
\[
A_n=\sum_{i,j}\binom ni^2\binom nj\binom ij\binom{n+i}n\binom{n+j}n 
\]
\[
A_n=\sum_{i,j}\binom ni^2\binom nj\binom{2n-j}n\binom{n+j}n\binom{2n-i}{n+j} 
\]
\[
A_n=\sum_{i,j}\binom ni^2\binom nj\binom{2n-j}n\binom{n+j}n\binom{n+i}{n+j} 
\]
\[
A_n=\sum_{i,j}\binom ni^2\binom nj^2\binom{2n-j}n\binom{n+i-j}i 
\]
\[
A_n=\sum_{i,j}\binom ni^2\binom nj\binom{n+j}n\binom{i+j}n\binom{n+i}{n-j} 
\]
\[
A_n=\sum_{i,j}\binom ni^2\binom nj\binom{i+j}j\binom{n+i-j}{n-j}\binom{n+i}{%
n-j} 
\]
\[
A_n=\sum_{i,j}\binom ni^2\binom nj\binom{i+j}j\binom{2i}n\binom{2i}{n-j} 
\]
\[
A_n=\sum_{i,j}\binom ni^2\binom nj\binom{2j}j\binom{2j}n\binom{i+j}{n-j} 
\]
\[
A_n=\sum_{i,j}\binom ni^2\binom nj\binom{i+j}j^2\binom{n+i}{n-j} 
\]
\[
A_n=\sum_{i,j}\binom ni^2\binom nj\binom{2n-i}n\binom{n+j}n\binom{n-i}j 
\]
\[
A_n=\sum_{i,j}\binom ni^2\binom nj\binom ij\binom{n+i}n\binom{2n-j}n 
\]
\[
A_n=\sum_{i,j}\binom ni^2\binom nj^2\binom{i+j}j\binom{3n-i-j}{n-j} 
\]
\[
A_n=\sum_{i,j}\binom ni^2\binom nj^2\binom{i+j}j\binom{3n-i-j}{2n-j} 
\]
\[
A_n=\sum_{i,j}(-1)^j\binom ni^2\binom{n+j}n^2\binom{i+j}n\binom{3n+1}{n-j} 
\]
\[
A_n=\sum_{i,j}(-1)^j\binom ni^2\binom{n+j}n^2\binom{n+i+j}n\binom{3n+1}{n-j} 
\]
\[
A_n=\sum_{i,j}(-1)^{i+j}\binom{2n-i}n^3\binom{n+j}n^2\binom{3n+1}{i-j} 
\]
\[
A_n=\sum_{i,j}(-1)^{i+j}\binom{2n-i}n^3\binom{n+j}n^2\binom{3n+1}{i-j} 
\]
\[
A_n=\sum_{i,j}(-1)^{i+j}\binom{2n-i}n^2\binom{n+j}n^3\binom{3n+1}{i-j} 
\]
\[
A_n=\sum_{i,j}\binom ni^2\binom nj\binom{n+i}n\binom{n+j}n\binom i{n-j} 
\]
\[
A_n=\sum_{i,j}\binom ni^2\binom nj\binom{n+i}n\binom{2n-j}n\binom{n+i-j}n 
\]
\[
A_n=\sum_{i,j}\binom ni\binom nj\binom{n+j}n\binom{i+j}n\binom{2j}n\binom
n{i-j} 
\]
\[
A_n=\sum_{i,j}\binom ni\binom nj\binom{n+j}n\binom{2n-i}n\binom{2n-j}n\binom
n{i-j} 
\]
\[
A_n=\sum_{i,j}\binom ni\binom nj\binom{n+i}n\binom{n+j}n\binom{2n-j}n\binom
n{i+j} 
\]
\[
A_n=\sum_{i,j}\binom ni\binom nj\binom{2n-i}n\binom{2n-j}n\binom{i+j}j\binom
n{i+j} 
\]
\[
"A_n"=(-1)^{n+1}\left\{ \frac{(3n)!}{n!^3}\right\}
^2\sum_k(-1)^kk^{-5}(n-k)^2(k+\frac n2)\binom{2n+k}{3n}^2\binom{n+k}n^{-5} 
\]
\[
\]

\textbf{196.} 
\[
A_n=\sum_{i,j}\binom ni^2\binom nj^2\binom{i+j}j\binom{n+i-j}n 
\]
\[
A_n=\sum_{i,j}\binom ni^2\binom nj\binom ij\binom{i+j}j\binom{n+i-j}{n-j} 
\]
\[
A_n=\sum_{i,j}\binom ni^2\binom nj\binom{2j}j\binom{n+i-j}{n-j}\binom{i+j}{2j%
} 
\]
\[
A_n=\sum_{i,j}\binom ni^2\binom nj^2\binom{i+j}n\binom{n+i-j}i 
\]
\[
A_n=\sum_{i,j}\binom ni^2\binom nj^2\binom{2n-i-j}n\binom{n+i-j}i 
\]
\[
A_n=\sum_{i.j}\binom ni^2\binom nj\binom ij\binom{n+i-j}{n-j}\binom{2n-i-j}{%
n-j} 
\]
\[
A_n=\sum_{i,j}\binom ni^2\binom nj^2\binom{n+i-j}n\binom{2n-i-j}{n-j} 
\]
\[
\]

\textbf{197.} 
\[
A_n=\sum_{i,j}\binom ni^2\binom nj^2\binom ij\binom{i+j}j 
\]
\[
A_n=\sum_{i,j}\binom ni^2\binom nj\binom ij\binom{2i}n\binom{2j}n 
\]
\[
A_n=\sum_{i,j}\binom ni^2\binom nj\binom{i+j}j\binom{n+i-j}{n-j}\binom
n{i+j} 
\]
\[
A_n=\sum_{i,j}\binom ni^2\binom nj\binom ij\binom{2i}n\binom i{n-j} 
\]
\[
A_n=\sum_{i,j}\binom ni^2\binom nj^2\binom{2j}j\binom{i+j}{2j} 
\]
\[
\]

\textbf{198.} 
\[
A_n=\sum_{i,j}\binom ni^2\binom nj^2\binom{i+j}j\binom{2n-i}n 
\]
\[
A_n=\sum_{i,j}\binom ni^2\binom nj^2\binom{n+i}n\binom{i+j}n 
\]
\[
A_n=\sum_{i,j}\binom ni^2\binom nj\binom ij\binom{n+i}n\binom{n+j}n 
\]
\[
A_n=\sum_{i,j}\binom ni^2\binom nj\binom{n+j}n\binom{n+i-j}{n-j}\binom{n+i}{%
n+j} 
\]
\[
A_n=\sum_{i,j}\binom ni^2\binom nj\binom{n+j}n^2\binom{n+i}{n+j} 
\]
\[
A_n=\sum_{i,j}\binom ni^2\binom nj\binom ij\binom{n+i}n\binom{n+i-j}{n-j} 
\]
\[
A_n=\sum_{i,j}\binom ni^2\binom nj\binom ij\binom{2n-i}n\binom{2n-j}n 
\]
\[
A_n=\sum_{i,j}\binom ni^2\binom nj\binom{i+j}j\binom{i+j}n\binom{n+i}{n-j} 
\]
\[
A_n=\sum_{i,j}\binom ni^2\binom nj\binom{i+j}j\binom{n+i-j}n\binom{n+i}{n-j} 
\]
\[
A_n=\sum_{i,j}\binom ni^2\binom nj\binom{i+j}j\binom{i+j}n\binom{n+i}{n-j} 
\]
\[
A_n=\sum_{i,j}\binom ni^2\binom nj\binom{2n-j}n\binom{i+j}n\binom{n+i}{n-j} 
\]
\[
A_n=\sum_{i,j}\binom ni^2\binom nj^2\binom{n+i}n\binom{n+i-j}n 
\]
\[
A_n=\sum_{i,j}\binom ni^2\binom nj^2\binom{2n-i}n\binom{n+i-j}n 
\]
\[
A_n=\sum_{i,j}\binom ni^2\binom nj^2\binom{2n-i-j}{n-j}\binom{n+j}n 
\]
\[
A_n=\sum_{i,j}\binom ni^2\binom nj\binom ij\binom{2n-j}n\binom{n+i-j}{n-j} 
\]
\[
A_n=\sum_{i,j}\binom ni^2\binom nj\binom{i+j}j\binom{2n-i-j}{n-j}\binom{n+i}{%
n-j} 
\]
\[
A_n=\sum_{i,j}\binom ni\binom nj\binom{n+j}n\binom{n+i-j}n^2\binom n{i-j} 
\]
\[
A_n=\sum_{i,j}\binom ni\binom nj\binom{n+j}n^2\binom{n+i-j}n\binom n{i-j} 
\]
\[
A_n=\sum_{i,j}\binom ni\binom nj\binom{2n-i}n^2\binom{n+j}n\binom n{i-j} 
\]
\[
A_n=\sum_{i,j}\binom ni\binom nj\binom{2n-i}n\binom{n+j}n\binom{i+j}j\binom
n{i+j} 
\]
\[
A_n=\sum_{i,j}\binom ni^2\binom nj^2\binom{i+j}n\binom{n+j}n 
\]
\[
A_n=\sum_{i,j}\binom ni^2\binom nj\binom{n+j}n\binom{i+j}j\binom i{n-j} 
\]
\[
A_n=\sum_{i,j}\binom ni^2\binom nj^2\binom{n+j}n\binom{2n-i-j}{n-j} 
\]
\[
A_n=\sum_{i,j}\binom ni^2\binom nj\binom{n+j}j\binom{n+i}{n-j}\binom{2n-i-j}{%
n-i} 
\]
\[
A_n=\sum_{i,j}\binom ni^2\binom nj^2\binom{n+j}n\binom{n+i-j}i 
\]
\[
A_n=\sum_{i,j}\binom ni\binom nj\binom{2n-i}n\binom{n+j}n\binom{i+j}j\binom
n{i+j} 
\]
\[
A_n=\sum_{i,j}\binom ni\binom nj\binom{n+i}n^2\binom{n+j}n\binom n{i+j} 
\]
\[
A_n=\sum_{i,j}\binom ni^2\binom nj\binom{2n-i}n\binom{n+j}n\binom{n-i}j 
\]
\[
"A_n"=\left\{ \frac{(3n)!}{n!^3}\right\} ^2\sum_k(-1)^k(\frac
n2+k)(k-n)k^{-6}\binom{2n+k}{3n}\binom{n+k}n^{-6} 
\]
\[
"A"_n=\sum_k(n-2k)\binom nk^5\binom{n+k}n\binom{2n-k}n 
\]
\[
\]

\textbf{202.} 
\[
A_n=\sum_{i,j}\binom ni^2\binom nj^2\binom{i+j}j\binom{n+i-j}{n-j} 
\]
\[
A_n=\sum_{i,j}\binom ni^2\binom nj^2\binom{n+i-j}n\binom{n+i-j}{n-j} 
\]
\[
A_n=\sum_{i,j}\binom ni^2\binom nj^2\binom{i+j}n\binom{2i}n 
\]
\[
A_n=\sum_{i,j}\binom ni^2\binom nj\binom ij\binom{n+i-j}{n-j}^2 
\]
\[
A_n=\sum_{i,j}\binom ni^2\binom nj\binom ij\binom{n+i-j}{n-j}\binom{2i}n 
\]
\[
A_n=\sum_{i,j}\binom ni^2\binom nj^2\binom{n+i-j}n\binom{2i}n 
\]
\[
A_n=\sum_{i,j}\binom ni^2\binom nj^2\binom{n+i-j}{n-j}\binom{2n-i-j}{n-j} 
\]
\[
A_n=\sum_{i,j}\binom ni^2\binom nj^2\binom{i+j}j\binom{i+j}n 
\]
\[
A_n=\sum_{i,j}\binom ni^2\binom nj^2\binom{i+j}j\binom{n+i-j}i 
\]
\[
A_n=\sum_{i,j}\binom ni^2\binom nj^2\binom{n+i-j}n\binom{n+i-j}i 
\]
\[
A_n=\sum_{i,j}\binom ni^2\binom nj\binom ij\binom{n+j}n\binom{2i}n 
\]
\[
\]

\textbf{208.} 
\[
A_n=\binom{2n}n^2\sum_k\binom nk^2\binom{n+2k}n 
\]
\[
A_n=\binom{2n}n\sum_k\binom nk\binom{n+k}n\binom{2n}k\binom{2k}n 
\]
\[
A_n=\binom{2n}n\sum_{i,j}\binom ni\binom nj\binom{i+j}j\binom{i+j}n\binom{2n%
}{i+j} 
\]
\[
A_n=\binom{2n}n\sum_{i,j}\binom ni\binom nj\binom{i+j}n\binom{n+i}{n-j}%
\binom{2n}{i+j} 
\]
\[
A_n=\binom{2n}n\sum_{i,j}\binom ni\binom nj\binom{n+i}n\binom{2i}n\binom
n{i-j} 
\]
\[
A_n=\binom{2n}n^2\sum_{i,j}\binom ni\binom nj\binom{n+i}{n-j}\binom n{i+j} 
\]
\[
A_n=\binom{2n}n^2\sum_{i,j}\binom ni^2\binom nj\binom{2i}{n-j} 
\]
\[
A_n=\binom{2n}n^2\sum_{i,j}(-1)^{i+j}\binom ni\binom nj\binom{n+2j}n\binom
n{2j-i} 
\]
\[
\]

\textbf{209.} 
\[
A_n=\binom{2n}n\sum_k\binom nk^2\binom{n+k}n\binom{n+2k}n 
\]
\[
A_n=\binom{2n}n\sum_k(-1)^{n+k}\binom nk\binom{n+k}n^3 
\]
\[
A_n=\binom{2n}n^3\sum_k(-1)^k\binom nk^3\binom{2n-k}n\binom{2n}k^{-1} 
\]
\[
A_n=\binom{2n}n^2\sum_k(-1)^k\binom nk^2\binom{2n-k}n^2\binom{2n}k^{-1} 
\]
\[
A_n=\sum_{i,j}(-1)^{i+j}\binom ni\binom nj\binom{i+j}n^4 
\]
\[
A_n=\binom{2n}n\sum_{i,j}\binom ni^2\binom nj\binom{n+j}n\binom{n+i}{n-j} 
\]
\[
A_n=\sum_{i,j}\binom ni^2\binom nj\binom{i+j}j\binom{n+2i}n\binom{n+i}{n-j} 
\]
\[
A_n=\binom{2n}n\sum_{i,j}\binom ni^2\binom nj\binom{2n-j}n\binom{n+i}{n+j} 
\]
\[
A_n=\binom{2n}n\sum_{i,j}\binom ni^2\binom nj\binom ij\binom{n+2i}n 
\]
\[
A_n=\binom{2n}n\sum_{i,j}\binom ni\binom nj\binom{n+2j}n\binom ij\binom{2j}{%
n-j} 
\]
\[
A_n=\binom{2n}n\sum_{i,j}\binom ni\binom nj\binom{n+2j}n\binom ij\binom{2j}i 
\]
\[
A_n=\binom{2n}n\sum_{i,j}\binom ni^2\binom nj\binom{2n-j}n\binom{n+i}j 
\]
\[
A_n=\binom{2n}n\sum_{i,j}\binom ni^2\binom nj\binom{2n-j}n\binom{2n-i}{n+j} 
\]
\[
A_n=\binom{2n}n\sum_{i,j}\binom ni^2\binom nj^2\binom{2n+i-j}i 
\]
\[
A_n=\binom{2n}n\sum_{i,j}\binom ni^2\binom nj^2\binom{2n+i-j}{n-j} 
\]
\[
A_n=\sum_{i,j}\binom ni\binom nj\binom{n+j}n\binom{i+j}n^2\binom{2n}{i+j} 
\]
\[
A_n=\binom{2n}n\sum_{i,j}\binom ni\binom nj\binom{2n-i}n\binom{2n-j}n\binom
n{i-j} 
\]
\[
A_n=\binom{2n}n\sum_{i,j}\binom ni\binom nj\binom{2n-j}n\binom{n+i-j}n\binom
n{i-j} 
\]
\[
A_n=\sum_{i,j}\binom ni\binom nj\binom{n+i}n\binom{2n-j}n\binom ij\binom{2n}%
i 
\]
\[
A_{i,j}=\binom{2n}n\sum_{i,j}\binom ni^2\binom nj\binom{n+2j}n\binom i{n-j} 
\]
\[
A_n=\binom{2n}n\sum_{i,j}\binom ni^2\binom nj\binom ij\binom{n+2i}n 
\]
\[
A_n=\sum_{i,j}\binom ni\binom nj\binom{2n-j}n\binom{n+i}n\binom{i+j}n\binom{%
2n}{i+j} 
\]
\[
A_n=\binom{2n}n\sum_{i,j}(-1)^{i+j}\binom{2n-i}n\binom{2n-j}n^3\binom{2n+1}{%
i-j} 
\]
\[
A_n=\sum_{i,j}(-1)^{i+j}\binom{2n-i}n^4\binom{n+j}n\binom{3n+1}{i-j} 
\]
\[
A_n=\sum_{i,j}(-1)^{i+j}\binom{2n-i}n\binom{n+j}n^3\binom{3n-i}n\binom{3n+1}{%
i-j} 
\]
\[
A_n=\sum_{i,j}(-1)^{i+j}\binom{2n-i}n^3\binom{n+j}n\binom{2n+j}n\binom{3n+1}{%
i-j} 
\]
\[
A_n=\binom{2n}n\sum_{i,j}\binom ni^2\binom nj^2\binom{n+i+j}n 
\]
\[
A_n=\binom{2n}n\sum_{i,j}(-1)^{i+j}\binom ni\binom nj\binom{n+j}n\binom{n+i+j%
}n^2 
\]
\[
A_n=\binom{2n}n\sum_{i,j}(-1)^{i+j}\binom ni\binom nj\binom{n+j}n^2\binom{%
n+i+j}n 
\]
\[
A_n=\binom{2n}n\sum_{i,j}(-1)^{i+j}\binom ni\binom nj\binom{i+j}j\binom{n+j}n%
\binom{n+i+j}n 
\]
\[
A_n=\binom{2n}n\sum_{i,j}\binom ni\binom nj\binom{n+i}n\binom{2n-j}n\binom
n{i+j} 
\]
\[
A_n=\binom{2n}n\sum_{i,j}\binom ni^2\binom nj\binom{n+i}n\binom{2i}{n-j} 
\]
\[
"A_n"=\binom{2n}n\sum_k(-1)^kk^{-4}\binom nk^{-3}\binom{n+k}n^{-1} 
\]
\[
"A_n"=(2n+1)^{-1}\sum_k\binom{n+k}n^4\binom{2n+k+1}k^{-1} 
\]
\[
\]

\textbf{210.} 
\[
A_n=\binom{2n}n\sum_k(-1)^k\binom{2n}k^4 
\]
\[
A_n=\binom{2n}n\sum_k(-1)^k\binom{2n}{n+k}^4 
\]
\[
A_n=\binom{2n}n^2\sum_{i,j}\binom nk\binom{2k}k\binom{2n+k}{2k} 
\]
\[
A_n=\binom{2n}n^2\sum_{i,j}\binom nk\binom{2n}k\binom{2n+k}{2n} 
\]
\[
A_n=\binom{2n}n^2(-1)^n\sum_{i,j}\binom nk\binom{2n}k\binom{3n-k}n 
\]
\[
A_n=\binom{2n}n^2\sum_k(-1)^{n+k}\binom{n+k}n\binom{2n}k\binom{2n+k}{2n} 
\]
\[
A_n=\binom{2n}n^2\sum_k(-1)^k\binom{2n}{n-k}\binom{2n+k}{2n}^2 
\]
\[
A_n=\binom{2n}n\sum_k(-1)^k\binom nk\binom{n+k}n\binom{2n+k}n^2 
\]
\[
A_n=(-1)^n\binom{2n}n^2\sum_k\binom nk^2\binom{2n-k}n\binom{3n-k}n\binom{n+k}%
n^{-1} 
\]
\[
A_n=\binom{2n}n\sum_{i,j}\binom ni^2\binom nj^2\binom{2n+i+j}{n+j} 
\]
\[
A_n=\binom{2n}n^2\sum_{i,j}\binom ni\binom nj\binom{2n+j}n\binom n{i+j} 
\]
\[
A_n=\binom{2n}n\sum_{i,j}\binom ni\binom nj\binom{2n-i}n\binom{n+j}n\binom{2n%
}{i+j} 
\]
\[
A_n=(-1)^n\binom{2n}n^2\sum_{i,j}\binom ni\binom nj\binom{i+j}j\binom{2n}{i+j%
} 
\]
\[
A_n=(-1)^n\binom{2n}n\sum_{i,j}\binom ni\binom nj\binom{i+j}n\binom{n+i+j}n%
\binom{2n}{i+j} 
\]
\[
A_n=\binom{2n}n^2\sum_{i,j}\binom ni^2\binom nj\binom{n+2i}{n+j} 
\]
\[
A_n=\binom{2n}n\sum_{i,j}\binom ni^2\binom nj\binom{2n-j}n\binom{2n+i}{n+j} 
\]
\[
A_n=\sum_{i,j}\binom ni^2\binom nj\binom{i+j}j\binom{3n-j}n\binom{2n+i}{n+j} 
\]
\[
A_n=\sum_{i,j}(-1)^{i+j}\binom{i+j}j\binom{2n-i}n\binom{n+j}n\binom{3n-i}n%
\binom{2n+j}n\binom n{i+j} 
\]
\[
A_n=\sum_{i,j}(-1)^{i+j}\binom{2n-i}n\binom{n+j}n\binom{2n+j}n^2\binom{3n-i}n%
\binom{3n+1}{i-j} 
\]
\[
A_n=\sum_{i,j}(-1)^{i+j}\binom{2n-i}n\binom{n+j}n\binom{2n+j}n\binom{3n-i}n^2%
\binom{3n+1}{i-j} 
\]
\[
A_n=\sum_{i,j}(-1)^{i+j}\binom{2n-i}n^2\binom{n+j}n\binom{3n-i}n^2\binom{3n+1%
}{i-j} 
\]
\[
A_n=\sum_{i,j}(-1)^{i+j}\binom{2n-i}n\binom{n+j}n^2\binom{2n+j}n^2\binom{3n+1%
}{i-j} 
\]
\[
A_n=\binom{2n}n\sum_{i,j}(-1)^{i+j}\binom{n+j}n\binom{2n+i}n\binom{2n+j}n^2%
\binom{3n+1}{i-j} 
\]
\[
A_n=\binom{2n}n\sum_{i,j}(-1)^{i+j}\binom ni\binom nj\binom{2n}{i+j}^3 
\]
\[
"A_n"=\binom{2n}n^3\sum_k(-1)^kk^{-3}\binom{n+k}n^{-2}\binom{2n+k}n^{-1}%
\binom{2n}k^{-1} 
\]
\[
"A_n"=\binom{2n}n^3\sum_kk^{-3}\binom nk^{-1}\binom{n+k}n^{-1}\binom{2n}%
k^{-1}\binom{2n+k}n^{-1} 
\]
\[
"A_n"=\sum_kk^{-2}(k-n)\binom kn^2\binom{2n-k}n^2\binom{n+k}n^{-1} 
\]
\[
"A"_n=\binom{2n}n^2\sum_k(n-2k)\binom nk^3\binom{2n+k}n\binom{3n-k}n 
\]
\[
\]

\textbf{211.} 
\[
A_n=\binom{2n}n^4\sum_{k=0}^n(-1)^k\binom nk^2\binom{2k}k\binom{4n-2k}{2n-k}%
\binom{n+k}n^{-2}\binom{2n}k^{-1} 
\]
\[
+\binom{2n}n^4\sum_{k=1}^n\binom nk^2\binom{2n+k}{2n}\binom{4n+2k}{2n+k}%
\binom{n+k}n^{-2}\binom{2k}k^{-1} 
\]
\[
A_n=\binom{2n}n^3\sum_k(-1)^{n+k}\binom nk\binom{2k}k\binom{2n-k}n\binom{%
4n-2k}{2n-k}\binom{n+k}n^{-2} 
\]
\[
"A_n"=256^n\binom{2n}n^3\sum_k(-1)^k\binom{2n+k}{2n}^6\binom{n+k}n^{-2}%
\binom{2n+k}n^{-2} 
\]
\[
"A_n"=256^{-n}\binom{2n}n\sum_k(-1)^k\binom{2n+k}{2n}^4 
\]
\[
"A_n"=256^n\binom{2n}n\sum_k(-1)^{n+k}k^{-4}\binom{2n}k^{-4} 
\]
\[
"A"_n=256^{-n}\binom{2n}n\sum_k(-1)^{n+k}\binom{n+k}{2n}^4 
\]
\[
"A_n"=256^{-n}\binom{2n}n\sum_k(-1)^k\binom{2n+k}n^4 
\]
\[
\]

\textbf{212.} 
\[
A_n=\sum_{i,j}\binom ni^2\binom nj^2\binom{i+j}j\binom{2n-i-j}n 
\]
\[
A_n=\sum_{i,j}\binom ni^2\binom nj^2\binom{i+j}n\binom{n+i-j}n 
\]
\[
A_n=\sum_{i,j}\binom ni^2\binom nj^2\binom{i+j}n\binom{2n-i-j}{n-j} 
\]
\[
A_n=\sum_{i,j}\binom ni^2\binom nj\binom ij\binom{i+j}n\binom{n+i-j}{n-j} 
\]
\[
A_n=\sum_{i,j}\binom ni^2\binom nj\binom ij\binom{i+j}n\binom{n+i-j}i 
\]
\[
A_n=\sum_{i,j}\binom ni^2\binom nj\binom ij\binom{n+i-j}i\binom{n+j-i}j 
\]
\[
A_n=\sum_{i,j}\binom ni^2\binom nj\binom ij\binom{n+i-j}{n-j}\binom{2n-i-j}n 
\]
\[
A_n=\sum_{i,j}\binom ni^2\binom nj\binom{i+j}j\binom i{n-j}\binom{2n-i-j}{n-i%
} 
\]
\[
A_n=\sum_{i,j}\binom ni^2\binom nj\binom{i+j}j\binom{2n-i-j}{n-j}\binom{n-i}%
j 
\]
\[
\,A_n=\sum_{i,j}(-1)^j4^{n-j}\binom ni^2\binom nj\binom{i+j}j\binom{2n-j}n%
\binom{2j}j 
\]
\[
A_n=\sum_{i,j}\binom ni^2\binom nj^2\binom{i+j}j\binom{2i}n 
\]
\[
A_n=\sum_{i,j}\binom ni^2\binom nj\binom ij\binom{i+j}j\binom{2j}n 
\]
\[
"A"_n=\sum_k(n-2k)\binom nk^5\binom{2k}k\binom{2n-2k}{n-k} 
\]
\[
\]

\textbf{213.} 
\[
A_n=\sum_{i,j}\binom ni^2\binom nj^2\binom{i+j}j\binom{2i}n 
\]
\[
A_n=\sum_{i,j}\binom ni^2\binom nj\binom ij\binom{2n-j}n\binom{2i}n 
\]
\[
A_n=\sum_{i,j}\binom ni\binom nj\binom{n+i}n\binom{n+j}n\binom{2n-i-j}{n-j}%
\binom n{i+j} 
\]
\[
A_n=\sum_{i,j}\binom ni\binom nj\binom{2n-i}n\binom{n+j}n\binom{2i}n\binom
n{i-j} 
\]
\[
A_n=\sum_{i,j}\binom ni\binom nj\binom{2n-i}n\binom{n+j}n\binom{n+i-j}{n-j}%
\binom n{i-j} 
\]
\[
A_n=\sum_{i,j}\binom ni^2\binom nj\binom{i+j}n\binom{2j}n\binom{n+i}{n-j} 
\]
\[
A_n=\sum_{i,j}\binom ni^2\binom nj^2\binom{n+i-j}n\binom{2i}n 
\]
\[
\]

\textbf{214.} 
\[
A_n=\binom{2n}n\sum_k\binom nk^2\binom{n+k}n\binom{3k}n 
\]
\[
A_n=\binom{2n}n\sum_{i,j}(-1)^{i+j}\binom ni\binom nj\binom{i+j}j^2\binom{i+j%
}n 
\]
\[
A_n=\binom{2n}n\sum_{i,j}(-1)^j\binom ni^2\binom{n+i}n\binom{3i+1}{n-j} 
\]
\[
\]

\textbf{215.} 
\[
A_n=\binom{2n}n^2\sum_k\binom nk^2\binom{4k}{2n} 
\]
\[
A_n=\binom{2n}n\sum_k(-1)^k4^{n-k}\binom nk\binom{n+k}n\binom{2k}n\binom{n+2k%
}k 
\]
\[
A_n=\binom{2n}n^2\sum_{i,j}(-1)^{i+j}\binom ni\binom nj\binom{4j}{2n}\binom
n{2j-i} 
\]
\[
\]

\textbf{217.} 
\[
A_n=\binom{2n}n\sum_k\binom nk^2\binom{3k}n\binom{2n-2k}n 
\]
\[
A_n=\binom{2n}n\sum_{i,j}(-1)^{n+i+j}\binom ni\binom nj\binom{i+j}n^2\binom{%
2n-i-j}n 
\]
\[
\]

\textbf{218.} 
\[
A_n=\binom{2n}n\sum_k\binom nk^2\binom{2k}n\binom{3k}n 
\]
\[
A_n=\binom{2n}n\sum_{i,j}\binom ni\binom nj\binom{n+j}n\binom{i+j}n\binom
n{i-j} 
\]
\[
A_n=\binom{2n}n\sum_{i,j}(-1)^{i+j}\binom ni\binom nj\binom{2n-i}n\binom{%
2n-i-j}n^2 
\]
\[
A_n=\binom{2n}n\sum_{i,j}(-1)^{i+j}\binom ni\binom nj\binom{n+j}n\binom{i+j}%
j^2 
\]
\[
\]

\textbf{219.} 
\[
A_n=\binom{2n}n^2\sum_k\binom nk^2\binom{3k}n 
\]
\[
A_n=\binom{2n}n^2\sum_{i,j}(-1)^{i+j}\binom ni\binom nj\binom{3j}n\binom{2n}{%
2j-i} 
\]
\[
\]

\textbf{220.} 
\[
A_n=\binom{2n}n^3\left\{ (-1)^n\sum_{k=0}^{[n/2]}\binom nk^2\binom n{2k}%
\binom{2n}{4k}^{-1}+\sum_{k=[n/2]+1}^n\binom nk^2\binom{4k}{2n}\binom{2k}%
n^{-1}\right\} 
\]
\[
A_n=\binom{2n}n^2\sum_{i,j}\binom ni\binom nj\binom{n+i}{n-j}\binom{2n}{i+j} 
\]
\[
\]

\textbf{222.} 
\[
A_n=\binom{2n}n^2\sum_k\binom nk\binom{2k}n\binom{2n}{n-k} 
\]
\[
A_n=\binom{2n}n\sum_k\binom{2n-k}k\binom{2n-2k}n\binom{2n}k^2 
\]
\[
A_n=\frac{(3n)!}{n!^3}\sum_k\binom nk\binom{2n}k\binom{2n}{n+k} 
\]
\[
A_n=\sum_{i,j}\binom ni^2\binom nj\binom{n+i-j}{n-j}\binom{2j}n\binom{2n+i}{%
n+j} 
\]
\[
A_n=\binom{2n}n^2\sum_{i,j}(-1)^{i+j}\binom ni^2\binom nj\binom n{2j-i} 
\]
\[
\]

\textbf{229.} 
\[
A_n=\sum_{i+j+k=n}\left( \frac{(2n)!}{i!j!k!}\right) ^2 
\]
\[
A_n=\sum_k\binom{2n}k^2\binom{4n-2k}{2n-k} 
\]
\[
\]

\textbf{232.} 
\[
A_n=\binom{2n}n^2\sum_k\binom nk^2\binom{3n}{n+k} 
\]
\[
A_n=\binom{2n}n^2\sum_k(-1)^{n+k}4^{n-k}\binom nk\binom{2n+k}n\binom{3n+2k}{%
n+k} 
\]
\[
A_n=\frac{(3n)!}{n!^3}\sum_k\binom nk\binom{2n}k\binom{2k}n 
\]
\[
A_n=\frac{(3n)!}{n!^3}\sum_k\binom nk\binom{2n}k\binom{2n}{n-k} 
\]
\[
A_n=32^{-n}(-1)^n\binom{2n}n^2\sum_k\binom nk\binom{n+k}n\binom{2n-k}n\binom{%
2n+2k}{n+k}\binom{4n-2k}{2n-k}\binom{2n}{2k}^{-2} 
\]
\[
A_n=\sum_{i,j}\binom ni\binom nj\binom{i+j}n^2\binom{3n}{i+j}\binom{2n}{i+j} 
\]
\[
A_n=\binom{2n}n^2\sum_{i,j}\binom ni^2\binom nj\binom{2n}{i+j} 
\]
\[
\]

\textbf{233.} 
\[
A_n=\binom{2n}n^3\sum_k\binom nk^2\binom{3n}{n+k}\binom{2n}{2k}^{-1} 
\]
\[
A_n=2^{-n}\binom{2n}n\sum_k\binom nk\binom{n+k}n\binom{2n-k}n\binom{2n+2k}{%
n+k}\binom{4n-2k}{2n-k}\binom{2n}{2k}^{-1} 
\]
\[
A_n=\binom{2n}n^2\sum_k(-1)^{n+k}4^{n-k}\binom nk\binom{2k}k\binom{n+k}n^2 
\]
\[
\]

\textbf{235.} 
\[
A_n=\sum_k\binom nk\binom{2k}k\binom{2n-2k}{n-k}\binom{2k}{n-k}\binom{2n-2k}%
k 
\]
\[
A_n=\sum_{i,j}(-1)^{n+j}4^{n-j}\binom ni^2\binom nj\binom{2i}n\binom{2j}j%
\binom{i+j}n 
\]
\[
\]

\textbf{237.} 
\[
A_n=8^{-n}\binom{2n}n^4\sum_k\binom nk^2\binom{3n}{n+k}\binom{2n}{2k}^{-2} 
\]
\[
A_n=\binom{2n}n^2\binom{3n}n\sum_k\binom{2n}k^2\binom{2n}{n-k}^2\binom{2n}{2k%
}^{-2}\binom{3n}{n+k}^{-1} 
\]
\[
A_n=\binom{2n}n^2\sum_k\binom nk^3\binom{2n+2k}{n+k}\binom{4n-2k}{2n-k}%
\binom{2n}k^{-1}\binom{2n}{n-k}^{-1} 
\]
\[
A_n=\binom{2n}n\sum_k(-1)^{n+k}4^{n-k}\binom nk\binom{n+k}n\binom{2k}n\binom{%
2n+2k}{n+k} 
\]
\[
A_n=2^{-n}\binom{2n}n^2\sum_k(-1)^k4^{n-k}\binom nk^2\binom{2n-2k}n\binom{2n%
}{n-k}\binom{2n}{2k}^{-1}\binom{2n-2k}{n-k}^{-1} 
\]
\[
\]

\textbf{239.} 
\[
A_n=\frac{(3n)!}{n!^3}\sum_k\binom nk\binom{2n+2k}{n+k}\binom{4n-2k}{2n-k} 
\]
\[
A_n=8^{-n}\binom{2n}n^2\binom{3n}n\sum_k\binom nk\binom{2n+2k}{n+k}\binom{%
4n-2k}{2n-k}\binom{2n}{2k}^{-1} 
\]
\[
\]

\textbf{240.} 
\[
A_n=\sum_k\binom nk\binom{2k}k\binom{2n-2k}{n-k}\binom{n+2k}n\binom{3n-2k}n 
\]
\[
A_n=\binom{2n}n^2\left\{ \sum_{k=0}^{[n/2]}\binom nk^3\binom{3n-2k}{2n}%
\binom{2n}{n+2k}^{-1}+\sum_{k=[n/2]+1}^n\binom nk^3\binom{2n}{n+2k}\binom{2n%
}{3n-2k}^{-1}\right\} 
\]
\[
\]

\textbf{241.} 
\[
A_n=\sum_k\binom{2n+2k}{n+k}\binom{4n-2k}{2n-k}\binom{n+k}{n-k}\binom{2n-k}k 
\]
\[
A_n=\binom{2n}n^2\sum_k(-1)^k\binom nk\binom{2n}k^2\binom{2n}{2k}^{-1} 
\]
\[
A_n=\binom{2n}n^2\sum_k(-1)^k\binom{2n-k}n\binom{2n}k^3\binom{2n}{2k}^{-1} 
\]
\[
A_n=\binom{2n}n^4\sum_k(-1)^k\binom nk\binom{2k}k\binom{2n-2k}{n-k}\binom{n+k%
}n^{-2} 
\]
\[
A_n=\binom{2n}n^3\sum_k(-1)^k\binom nk\binom{2n}k^2\binom{2n}{2k}^{-1} 
\]
\[
\]

\textbf{243.} 
\[
A_n=\sum_{i,j}\binom ni\binom nj\binom{n+i}n\binom{n+j}n\binom{n+i+j}n\binom
n{i-j} 
\]
\[
A_n=\sum_{i,j}\binom ni^2\binom nj\binom{2n-j}n\binom{n+i-j}{n-j}\binom{n+2i%
}{n+j} 
\]
\[
\]

\textbf{248.} 
\[
A_n=(-1)^n\sum_k\binom nk^3\binom{n+k}n\binom{2n-k}n\binom{2k}k\binom{2n-2k}{%
n-k} 
\]
\[
A_n=\sum_{i,j}(-1)^j4^{n-j}\binom ni^2\binom nj\binom{i+j}j\binom{n+i}n%
\binom{2j}j 
\]
\[
A_n=\sum_{i,j}(-1)^j4^{n-j}\binom ni^2\binom nj\binom{i+j}n\binom{n+i}n%
\binom{2j}j 
\]
\[
\]

\textbf{256.} 
\[
A_n=\binom{2n}n\sum_k(-1)^k4^{n-k}\binom nk\binom{n+k}n^2\binom{2n+k}{n+k} 
\]
\[
A_n=\sum_k\binom nk^3\binom{n+k}n\binom{2n-k}n\binom{2n+2k}{n+k}\binom{4n-2k%
}{2n-k}\cdot 
\]
\[
\left\{ 1+k(-4H_k+4H_{n-k}-H_{n+k}+H_{2n-k}+2H_{2n+2k}-2H_{4n-2k})\right\} 
\]
\[
A_n=\binom{2n}n^2\sum_k\binom nk\binom{2n}k\binom{2n}{n-k}\binom{2n-k}k%
\binom{n+k}{n-k}\cdot 
\]
\[
\left\{ 1+k(-3H_k+3H_{n-k}-2H_{2k}+2H_{2n-2k})\right\} 
\]
\[
A_n=\binom{2n}n^2\sum_{i,j}\binom ni\binom nj\binom{n+i}{n-j}\binom{2n}{i+j} 
\]
\[
\]

\textbf{258.} 
\[
A_n=\binom{2n}n^2\sum_k4^{n-k}\binom nk^2\binom{n+k}n\binom{2n+2k}{n+k}%
\binom{2n}{2k}^{-1} 
\]
\[
A_n=\sum_k\binom nk\binom{n+k}n\binom{2n-k}n\binom{2k}k\binom{2n-2k}{n-k}%
\binom{2n+2k}{n+k}\binom{4n-2k}{2n-k} 
\]
\[
(1+k(-4H_k+4H_{n-k}-H_{n+k}+H_{2n-k}+2H_{2k}-2H_{2n-2k}+2H_{2n+2k}-2H_{4n-2k})) 
\]
\[
A_n=\binom{2n}n^2\sum_k\binom{2k}k\binom{2n}k\binom{4n-2k}{2n-k} 
\]
\[
\]

\textbf{277.} 
\[
A_n=\binom{2n}n^3(-1)^n\sum_k4^{n-k}\binom nk\binom{2n}{n+k}\binom{2n+2k}{n+k%
}\binom{2n}{2k}^{-1} 
\]
\[
A_n=\binom{2n}n^2\sum_k\binom nk\binom{2k}k\binom{2n-2k}{n-k}\binom{2n+2k}{%
n+k}\binom{4n-2k}{2n-k} 
\]
\[
(1+k(-3H_k+3H_{n-k}+2H_{2k}-2H_{2n-2k}+2H_{2n-k}-2H_{n+k}+2H_{2n+2k}-2H_{4n-2k})) 
\]
\[
\]

\textbf{279.} 
\[
A_n=3\sum_{k=[2(n+1)/3]}^n(-1)^k\binom nk\binom n{3n-3k}\binom{3k}n^{-1}%
\frac{n-2k}{n-3k}\frac{(3k)!}{k!^3}\frac{(3n-3k)!}{(n-k)!^3} 
\]
\[
A_n=\sum_{i,j}(-1)^{n+k}3^{n-3i}\binom n{3i}\binom nj\binom ij\binom{n+j}n%
\frac{(3i)!}{i!^3} 
\]
\[
A_n=\sum_{i,j}(-1)^{n+k}3^{n-3i}\binom n{3i}\binom nj^2\binom{i+j}n\frac{%
(3i)!}{i!^3} 
\]
\[
\]

\textbf{284.} 
\[
A_n=\sum_{i,j}\binom ni\binom nj\binom{n+i}n\binom{n+j}n\binom{n+i-j}{n-j}%
\binom n{i+j} 
\]
\[
A_n=\sum_{i,j}\binom ni\binom nj\binom{2n-i}n\binom{n+j}n\binom{i+j}j\binom
n{i-j} 
\]
\[
\]

\textbf{287.} 
\[
A_n=\binom{2n}n\sum_{i,j}\binom ni\binom nj\binom{n+j}n\binom{i+j}n\binom
n{i-j} 
\]
\[
A_n=\sum_{i,j}\binom ni\binom nj\binom{i+j}n^2\binom{n+i}{n-j}\binom{2n}{i+j}
\]
\[
A_n=\sum_{i,j}\binom ni\binom nj\binom{i+j}j\binom{i+j}n\binom{2n-i}n\binom{%
2n}{i+j} 
\]
\[
A_n=\binom{2n}n\sum_{i,j}\binom ni\binom nj\binom{n+j}n\binom{n+i-j}{n-j}%
\binom n{i-j} 
\]
\[
A_n=\binom{2n}n\sum_{i,j}\binom ni\binom nj\binom{n+j}n\binom n{i+j}\binom{%
i+2j}n 
\]
\[
A_n=\binom{2n}n\sum_{i,j}\binom ni^2\binom nj\binom{i+j}j\binom{2i}{n-j} 
\]
\[
A_n=\binom{2n}n\sum_{i,j}\binom ni\binom nj\binom ij\binom{n+i}n\binom{2j}i 
\]
\[
A_n=\binom{2n}n\sum_{i,j}\binom ni^2\binom nj\binom{2n-j}n\binom{2i}{n+j} 
\]
\[
A_n=\binom{2n}n\sum_{i,j}\binom ni^2\binom nj\binom{2j}j\binom{i+j}{n-j} 
\]
\[
A_n=\binom{2n}n\sum_{i,j}\binom ni^2\binom nj^2\binom{2i+j}{n+j} 
\]
\[
\]

\textbf{288.} 
\[
A_n=\binom{2n}n^2\sum_k4^{n-k}\binom nk^2\binom{2n+k}n\binom{4n+2k}{2n+k}%
\binom{2n}{2k}^{-1} 
\]
\[
A_n=\binom{2n}n\binom{3n}n\binom{4n}n\frac{(4k)!}{(2k)!k!^2}\frac{(4n-4k)!}{%
(2n-2k)!(n-k)!^2}\binom{n+k}n^{-1}\binom{2n-k}n^{-1} 
\]
\[
\]

\textbf{292.} 
\[
A_n=\binom{2n}n\sum_k(-1)^{n+k}4^{n-k}\binom nk\binom{2n+k}n^2\binom{4n+2k}{%
2n+k} 
\]
\[
A_n=\binom{2n}n^3\sum_k4^{2n-k}\binom nk\binom{n+k}n\binom{2n-2k}{n-k}^{-1} 
\]
\[
\]

\textbf{293.} 
\[
A_n=\binom{2n}n\sum_k(-1)^{n+k}4^{n-k}\binom nk\binom{2k}k\binom{n+k}n^2 
\]
\[
A_n=\sum_{i,j}\binom ni\binom nj\binom{i+j}n\binom{2i}n\binom{2j}n\binom{2n}{%
i+j} 
\]
\[
\]

\textbf{301.} 
\[
A_n=\sum_{i,j}\binom ni\binom nj\binom{n+j}n\binom{i+j}j\binom{2i}n\binom
n{i-j} 
\]
\[
A_n=\sum_{i,j}\binom ni\binom nj\binom{n+j}n\binom{i+j}n\binom{2i}n\binom
n{i-j} 
\]
\[
\]

\textbf{303.} 
\[
A_n=\sum_{i,j}\binom ni^2\binom nj\binom{n+i}n\binom{2j}j\binom{i+j}{n-j} 
\]
\[
A_n=\sum_{i,j}\binom ni^2\binom nj\binom{i+j}j\binom{n+j}n\binom{2i}{n-j} 
\]
\[
A_n=\sum_{i,j}\binom ni^2\binom nj\binom{n+j}n\binom{2i}n\binom{2i}{n-j} 
\]
\[
A_n=\sum_{i,j}\binom ni\binom nj\binom{i+j}n\binom{2i}n\binom{2j}n\binom
n{i-j} 
\]
\[
A_n=\sum_{i,j}\binom ni^2\binom nj^2\binom{n+j}n\binom{2i+j}{n+j} 
\]
\[
A_n=\sum_{i,j}\binom ni^2\binom nj^2\binom{2i}n\binom{2i+j}j 
\]
\[
\]

\textbf{306.} 
\[
A_n=\sum_{i,j}\binom ni\binom nj\binom{n+i}n\binom{n+j}n\binom{2i}n\binom
n{i-j} 
\]
\[
A_n=\sum_{i,j}\binom ni^2\binom{n+j}n^2\binom{n+i+j}n\binom{3n+1}{n-2j} 
\]
\[
\]

\textbf{307.} 
\[
A_n=\binom{2n}n\sum_{i,j}\binom ni^2\binom nj^2\binom{2i}{n+j} 
\]
\[
A_n=\binom{2n}n\sum_{i,j}\binom ni^2\binom nj^2\binom{2i}{n-j} 
\]
\[
A_n=\binom{2n}n\sum_{i,j}\binom ni^2\binom nj\binom ij\binom n{i-j} 
\]
\[
A_n=\binom{2n}n\sum_{i,j}\binom ni^2\binom nj\binom ij\binom{2j}i 
\]

\textbf{308.} 
\[
A_n=\sum_{i,j}\binom ni\binom ij\binom{n+i-j}n\binom{2j}j\binom{2j}{i-j} 
\]
\[
A_n=\sum_{i,j}\binom ni\binom nj\binom{n+i-j}{n-j}\binom{2j}j\binom j{i-j} 
\]
\[
A_n=\binom{2n}n\sum_{i,j}\binom ni^2\binom nj\binom{2n-i-j}n\binom{2i}{n-j} 
\]
\[
A_n=\binom{2n}n\sum_{i,j}\binom ni\binom nj\binom{i+j}j\binom n{i+j}\binom{%
i+2j}n 
\]
\[
\]

\textbf{309.} 
\[
A_n=\binom{2n}n\sum_{i,j}\binom ni^2\binom nj\binom{n+j}n\binom{2i}{n-j} 
\]
\[
A_n=\binom{2n}n\sum_{i,j}\binom ni^2\binom nj^2\binom{2i+j}j 
\]
\[
A_n=\binom{2n}n\sum_{i,j}(-1)^{i+j}\binom ni\binom nj\binom{i+j}j\binom{i+j}%
n^2 
\]
\[
A_n=\sum_{i,j}\binom ni\binom nj\binom{i+j}j\binom{i+j}n^2\binom{2n}{i+j} 
\]
\[
\]

\textbf{310.} 
\[
A_n=\sum_{i,j}\binom ni^2\binom nj\binom{n+i}n\binom{n+j}n\binom{2i}{n-j} 
\]
\[
A_n=\sum_{i,j}\binom ni^2\binom nj^2\binom{n+i}n\binom{2i+j}j 
\]
\[
\]

\textbf{312. } 
\[
A_n=\sum_{i,j}\binom ni^2\binom nj\binom{n+i-j}{n-j}\binom{2j}n\binom{2i}{n+j%
} 
\]
\[
A_n=\sum_{i,j}\binom ni^2\binom nj\binom{2n-i-j}{n-j}\binom{2j}n\binom{n+i}{%
n+j} 
\]
\[
\]

\textbf{316.} 
\[
A_n=\binom{2n}n^2\sum_{i,j}\binom ni^2\binom nj\binom{3i}{n+j} 
\]
\[
A_n=\binom{2n}n^2\sum_{i,j}(-1)^{i+j}\binom ni\binom nj\binom{n+3j}{2n}%
\binom n{2j-i} 
\]
\[
\]

\textbf{318.} 
\[
A_n=\sum_{i,j}\binom ni^2\binom nj^2\binom{n+i-j}{n-j}\binom{2i+j}{n+j} 
\]
\[
A_n=\sum_{i,j}\binom ni^2\binom nj\binom{i+j}j\binom{2j}j\binom{i+j}{n-j} 
\]
\[
A_n=\sum_{i,j}\binom ni^2\binom nj\binom{i+j}j^2\binom{2i}{n-j} 
\]
\[
\]

\textbf{332.} 
\[
A_n=2^{-n}\binom{2n}n\sum_k(-1)^{n+k}4^{n-k}\binom nk\binom{2k}k\binom{2n-2k}%
n\binom{n+k}n 
\]
\[
A_n=2^{-n}\binom{2n}n^2\sum_k(-1)^{n+k}4^{n-k}\binom nk^2\binom{n-k}k\binom{%
n+k}n\binom{2n}{2k}^{-1} 
\]
\[
\]

\textbf{Some other sums.} 
\[
\]

\textbf{(i) } 
\[
\sum_{i,j}\binom ni\binom nj\binom{i+j}n^2\binom{2n}{i+j}=\binom{2n}n^3 
\]
\[
\]

\textbf{(ii)} 
\[
\sum_{i,j}\binom ni\binom nj\binom{i+j}n^2\binom{n+i+j}{2n}\binom{2n}{i+j}=%
\binom{2n}n^4 
\]
\[
\]

\textbf{(iii)} 
\[
\sum_{i,j}\binom ni\binom nj\binom{i+j}n^2\binom{3n}{i+j}=\binom{2n}n\binom{%
3n}n^2 
\]
\[
\]

\textbf{(iv)} 
\[
\sum_{i,j}\binom ni\binom nj\binom{i+j}n^2\binom{2n+i+j}{2n}\binom{2n}{i+j}=%
\binom{2n}n^2\binom{3n}n^2 
\]
\[
\]

\textbf{(v)} 
\[
\sum_{i,j}(-1)^{i+j}\binom{2n-i}n\binom{n+j}n\binom{2n+j}n\binom{3n-i}n%
\binom{2n+1}{i-j}=\binom{2n}n^2\binom{3n}n 
\]
\[
\]

\textbf{(vi)} 
\[
\sum_{i,j}(-1)^{i+j}\binom{2n-i}n^2\binom{n+j}n\binom{2n+j}n\binom{3n-i}n%
\binom{2n+1}{i-j}=\binom{2n}n^3\binom{3n}n 
\]
\[
\]

\textbf{(vii)} 
\[
\sum_{i,j}(-1)^{i+j}\binom{2n-i}n^2\binom{n+j}n^2\binom{2n+j}n\binom{3n-i}n%
\binom{3n+1}{i-j}=\binom{2n}n^5 
\]
\[
\]

\textbf{(viii)} 
\[
\sum_{i,j}(-1)^{i+j}\binom{2n-i}n\binom{n+j}n\binom{2n+j}n\binom{3n-i}n%
\binom{3n+1}{i-j}=\binom{2n}n^3 
\]
\[
\]

\textbf{(ix)} 
\[
\sum_{i,j}\binom ni\binom nj\binom{3n}{i+j}=\binom{5n}{2n} 
\]
\[
\]

\textbf{(x)} 
\[
\sum_{i,j}\binom ni\binom nj\binom{i+j}n\binom{2n}{i+j}=\binom{2n}n\binom{3n}%
n 
\]
\[
\]

\textbf{(xi)} 
\[
\sum_{i,j}\binom ni\binom nj\binom{n+j}n\binom{2n}{i+j}=\binom{3n}n^2 
\]
\[
\]

\textbf{(xii)} 
\[
\sum_{i,j}(-1)^{i+j}\binom{2n-i}n\binom{n+j}n\binom{2n}{i+j}^2\binom{2n+1}{%
i-j}=(-1)^n\binom{2n}n^3 
\]
\[
\]

\textbf{(xii).} 
\[
\sum_{i,j}(-1)^{n+i+j}\binom{2n-i}n\binom{n+j}n\binom{2n}{i+j}\binom{3n+1}{%
i-j}=2^n\binom{2n}n 
\]
\[
\]

\textbf{(xiii)} 
\[
\sum_{i,j}(-1)^i\binom nj^2\binom{2n-i}n\binom{2n-j}n\binom{3n-i}n\binom{%
n+j-i}n\binom{3n+1}i=\binom{2n}n^4 
\]
\[
\]

\textbf{(xiv)} 
\[
\sum_{i,j}\binom ni^2\binom nj\binom{i+j}j\binom{2n}j=\binom{2n}n^3 
\]
\[
\]

\textbf{(xv)} 
\[
\sum_{i,j}(-1)^j\binom ni^2\binom nj\binom{n+i+j}n\binom{2n+1}{n-j}=(-1)^n 
\]
\[
\]

\textbf{(xvi)} 
\[
\sum_k\binom nk^2\binom{2n+k}n\binom{3n+k}n=\binom{2n}n\binom{3n}n^2 
\]
\[
\]

\textbf{(xvii)} 
\[
\sum_{i,j}\binom ni\binom nj\binom{i+j}n\binom{3n-i-j}n\binom{3n}{i+j}=%
\binom{2n}n\binom{3n}n^2 
\]
\[
\]

\textbf{(xviii)} 
\[
\sum_{i,j}\binom ni\binom nj\binom{i+j}n^2\binom{3n-i-j}n\binom{3n}{i+j}=%
\binom{2n}n^3\binom{3n}n 
\]
\[
\]

\textbf{(ixx)} 
\[
\sum_{i,j}\binom ni\binom nj\binom{i+j}j\binom{3n-i}n\binom n{i+j}=\binom{2n}%
n^3 
\]
\[
\]

\textbf{(xx)} 
\[
\sum_{i,j}(-1)^{i+j}\binom{2n-i}n\binom{n+j}n\binom{3n-i}n^2\binom{2n+j}n%
\binom{2n+1}{i-j}=\binom{2n}n^2\binom{3n}n^2 
\]
\[
\]

\textbf{(xxi)} 
\[
\binom{2n}n^2\sum_k\binom nk^4\binom{2n}{2k}^{-3}=4^n\sum_k\binom nk^4\binom{%
2n}{2k}^{-1} 
\]
\[
\]

\textbf{(xxii)} 
\[
\sum_{i,j}\binom ni\binom nj\binom{2n+j}n\binom{2j}j\binom{2j}{i-j}=\binom{2n%
}n^3 
\]
\[
\]

\textbf{(xxiii)} 
\[
\sum_{i,j}\binom ni^2\binom nj\binom{i+j}j\binom{2n+i}{2n-j}=\binom{2n}n^2%
\binom{3n}n 
\]
\[
\]

\textbf{(xxiv)} 
\[
\sum_{i,j}(-1)^{i+j+n}\binom{2n-i}n\binom{n+i}n\binom{n+j}n\binom{3n+1}{i-j}=%
\binom{2n}n 
\]
\[
\]

\textbf{(xxv)} 
\[
\sum_{i,j}\binom ni\binom nj\binom{2n-j}n\binom{i+j}j\binom{2n}i\binom{2n}j=%
\binom{2n}n^4 
\]
\[
\]

\textbf{(xxvi)} 
\[
\sum_{i,j}\binom ni\binom nj\binom{2n-j}n\binom{n+i-j}{n-j}\binom{2n}i\binom{%
2n}j=\binom{2n}n^4 
\]

\[
\]

\textbf{(xxxvii)} 
\[
\sum_{i,j}(-1)^{i+j}\binom ni\binom nj\binom{n+i}n\binom{2n-i}n\binom{2j}%
i=2^n\binom{2n}n 
\]
\[
\]

\textbf{(xxxviii)} 
\[
\sum_{i,j}(-1)^{i+j}\binom ni\binom nj\binom{n+i}n\binom{n+i}n\binom{2j}i=2^n%
\binom{2n}n 
\]
\[
\]

\textbf{(ixxxx)} 
\[
A_n=\sum_{i,j}\binom ni\binom nj\binom{2n+j}n\binom ij\binom{2j}i=\binom{2n}%
n^3 
\]
\[
\]

\textbf{(xxxx)} 
\[
\sum_{i,j}(-1)^{i+j}\binom ni\binom nj\binom ij\binom{3j}{2n}=\binom{3n}n 
\]
\[
\]

(xxxxi) 
\[
\sum_{i,j}(-1)^{i+j}\binom ni\binom nj^2\binom ij\binom{3j}{2n}=\binom{2n}n%
\binom{3n}n 
\]
\[
\]

\textbf{(xxxxii)} 
\[
\sum_{i,j}(-1)^{i+j}\binom ni\binom nj\binom{n+j}n^2\binom ij\binom{3j}{2n}=%
\binom{2n}n^2\binom{3n}n 
\]
\[
\]

\textbf{(xxxxiii)} 
\[
\sum_{i,j}(-1)^{i+j}\binom ni\binom nj\binom{n+j}n\binom{3i+1}j=3^n\binom{2n}%
n 
\]
\[
\]

\textbf{(xxxxiv)} 
\[
\sum_{i,j}(-1)^{i+j}\binom ni\binom nj\binom{n+i}n^2\binom{2j}i=2^n\binom{2n}%
n^2 
\]
\[
\]

\textbf{(xxxxv)} 
\[
\sum_{i,j}\binom ni^2\binom nj\binom{n+i}{n+j}=\binom{2n}n^2 
\]
\[
\]

\textbf{(xxxxvi)} 
\[
\sum_{i,j}(-1)^{i+j}\binom{2n-i}n^2\binom{n+j}n\binom{3n+1}{i-j}=\binom{2n}n 
\]
\[
\]

\textbf{(xxxxvii)} 
\[
\sum_{i,j}(-1)^{i+j}\binom ni\binom nj\binom{i+j}j\binom{3n-2j}n=(-2)^n 
\]
\[
\]

(\textbf{xxxxviii)} 
\[
\sum_{i,j}(-1)^{i+j}\binom{2n-i}n\binom{n+i}n\binom{2n+j}n\binom{3n-i}n%
\binom{3n+1}{i-j}=(-1)^n\binom{2n}n^3 
\]
\[
\]

\textbf{(xxxxviiii)} 
\[
\sum_{i,j}(-1)^{i+j}\binom{2n+i}n^2\binom{n+j}n\binom{2n+j}n\binom{2n+1}{i-j}%
=\binom{2n}n^3 
\]
\[
\]

\textbf{(L)} 
\[
\sum_{i,j}(-1)^{i+j}\binom{2n-i}n\binom{n+j}n\binom{2n}{i+j}\binom{2n+1}{i-j}%
=(-1)^n\binom{2n}n 
\]
\[
\]

\textbf{(Li)} 
\[
\sum_{i,j}\binom ni^2\binom nj\binom{n+i-j}{n-j}\binom{2n+i}{n+j}=\binom{2n}%
n^2\binom{3n}n 
\]
\[
\]

\textbf{(Lii)} 
\[
\sum_{i,j}\binom ni\binom nj\binom{i+j}j\binom{3n}{i+j}=\binom{3n}n\binom{4n}%
n 
\]
\[
\]

\textbf{(Liii)} 
\[
\sum_{i,j}\binom ni\binom nj\binom{i+j}j^2\binom{2n}{i+j}=\binom{2n}n^3 
\]
\[
\]

\textbf{(Liv)} 
\[
\sum_{i,j}\binom ni^2\binom nj\binom{2n+i}{n+j}=\binom{3n}n^2 
\]
\[
\]

\textbf{(Lv)} 
\[
\sum_{i,j}\binom ni^2\binom nj\binom{2n+i}{i+j}=\binom{3n}n^2 
\]
\[
\]

\textbf{(Lvi)} 
\[
\sum_{i,j}\binom ni^2\binom nj\binom{i+j}j\binom{2n+i}{i+j}=\binom{2n}n^2%
\binom{3n}n 
\]
\[
\]

\textbf{(Lvii)} 
\[
\sum_{i,j}\binom ni^2\binom nj\binom{n+i-j}{n-j}\binom{2n+i}{n+j}=\binom{2n}%
n^2\binom{3n}n 
\]
\[
\]

\textbf{(Lviii)} 
\[
\sum_{i,j}\binom ni^2\binom nj\binom{n+j}n\binom{n+i-j}{n-j}\binom{2n+i}{n+j}%
=\binom{2n}n^4 
\]
\[
\]

\textbf{(Lviiii)} 
\[
\sum_{i,j}\binom ni^2\binom nj\binom{2n-j}n\binom{i+j}j\binom{2n+i}{i+j}=%
\binom{2n}n^4 
\]
\[
\]

\textbf{(Lx)} 
\[
\sum_{i,j}\binom ni^2\binom nj\binom{2n+i-j}{n-j}\binom{3n+1}{n+j}=\binom{3n}%
n^3 
\]
\[
\]

\textbf{(Lxi)} 
\[
\sum_{i,j}\binom ni^2\binom nj\binom{2n-j}n\binom{n+i-j}{n-j}\binom{3n}{n+j}=%
\binom{2n}n^3\binom{3n}n 
\]
\[
\]

\textbf{(Lxii)} 
\[
\sum_k(-1)^{n+k}\binom{n+k}n^2\binom{2n}{n+k}^3\binom{2n}{2k}^{-2}=16^n%
\binom{2n}n 
\]
\[
\]

\textbf{(Lxiii)} 
\[
\sum_k(-1)^{n+k}\binom nk\binom{n+k}n^2=\sum_k\binom nk^2\binom{n+k}n 
\]
\[
\]

\textbf{(Lxiv)} 
\[
\sum_{i,j}(-1)^{i+j}\binom{i+j}j\binom{2n-j}n\binom{n+j}n\binom n{i+j}=(-1)^n%
\binom{2n}n 
\]
\[
\]

\textbf{(Lxv)} 
\[
\sum_{i,j}(-1)^{i+j}\binom ni\binom nj\binom{2n-j}n\binom{3n-j}n\binom
n{2j-i}=\binom{2n}n^3 
\]
\[
\]

\textbf{(Lxvi)} 
\[
\sum_{i,j}(-1)^{i+j}\binom ni\binom nj\binom{2n}{i+j}^2=(-1)^n\binom{2n}n%
\binom{3n}n 
\]
\[
\]

\textbf{(Lxvii)} 
\[
\sum_{i,j}\binom ni^2\binom nj\binom{3n-i}{n+j}=\binom{3n}n^2 
\]
\[
\]

\textbf{(Lxviii)} 
\[
\sum_{i,j}(-1)^{i+j}\binom{2n-i}n\binom{n+j}n\binom{2n+j}n\binom{2i}n\binom{%
4n+1}{i-j}=2^n\binom{2n}n^2 
\]
\[
\]

\textbf{(Lxix)} 
\[
\sum_{i,j}(-1)^{i+j}\binom{2n-i}n^2\binom{n+j}n\binom{2n+j}n\binom{3n-i}n%
\binom{4n+1}{i-j}=\binom{2n}n^4 
\]
\[
\]

\textbf{(Lxx)} 
\[
\sum_{i,j}(-1)^{i+j}\binom{2n-i}n\binom{n+j}n^2\binom{2n-j}n\binom{3n-i}n%
\binom{2n+1}{i-j}=\binom{2n}n^3 
\]
\[
\]

\textbf{(Lxxi)} 
\[
\sum_{i,j}(-1)^{i+j}\binom{2n-i}n\binom{n+j}n\binom{2n-j}n\binom{3n-i}n%
\binom{2n+j}n\binom{2n+1}{i-j}=\binom{2n}n^2\binom{3n}n 
\]
\[
\]

\textbf{(Lxxii)} 
\[
\sum_{i,j}(-1)^{n+i+j}\binom ni\binom nj\binom ij\binom{n+j}n\binom{n-2j}%
j=2^{-1}\binom{2n}n^2 
\]
\[
\]

\textbf{(Lxxiii)} 
\[
\sum_{i,j}(-1)^{n+i+j}\binom ni\binom nj\binom{n+i-j}i\binom{n+j}n\binom{n-2j%
}j=1 
\]
\[
\]

\textbf{(Lxxiv)} 
\[
\sum_{i,j}(-1)^{n+i+j}\binom ni\binom nj\binom ij\binom{2n-j}n\binom{n-2j}%
j=2^{-1}\binom{2n}n 
\]
\[
\]

\textbf{(Lxxv)} 
\[
\sum_{i,j}(-1)^{n+i+j}\binom ni\binom nj\binom ij\binom{2n+j}n\binom{n-2j}%
j=2^{-1}\binom{2n}n\binom{3n}n 
\]
\[
\]

\textbf{(Lxxvi)} 
\[
\sum_{i,j}(-1)^{n+i+j}\binom ni\binom nj\binom{n+i-j}i\binom{2n-j}n\binom{%
n-2j}j=\binom{2n}n 
\]
\[
\]

\textbf{(Lxxvii)} 
\[
\sum_{i,j}(-1)^{n+i+j}\binom ni\binom nj\binom ij\binom{3n-j}n\binom{n-2j}%
j=2^{-1}\binom{2n}n^2 
\]
\[
\]

\textbf{(Lxxviii)} 
\[
\sum_{i,j}(-1)^{n+i+j}\binom ni\binom nj\binom ij\binom{kn+j}n\binom{n-2j}%
j=2^{-1}\binom{2n}n\binom{(k+1)n}n 
\]
\[
\]

\textbf{(Lxxix)} 
\[
\sum_{i,j}(-1)^{n+i+j}\binom ni\binom nj\binom ij\binom{n+j}n^2\binom{n-2j}%
j=2^{-1}\binom{2n}n^3 
\]
\[
\]

\textbf{(Lxxx)} 
\[
\sum_{i,j}(-1)^{n+i+j}\binom ni\binom nj\binom ij\binom{n+i}{n-j}^2\binom{%
n-2j}j=2^{-1}\binom{2n}n 
\]
\[
\]

\textbf{(Lxxxi)} 
\[
\sum_{i,j}\binom ni^2\binom nj\binom{kn+i}{n+j}=\binom{(k+1)n}n^2 
\]
\[
\]

\textbf{(Lxxxii)} 
\[
\sum_{i,j}(-1)^{n+i+j}\binom{2n-i}n^2\binom{n+j}n\binom{n+i}{n-j}=\binom{2n}%
n^2 
\]
\[
\]

(\textbf{Lxxxiii)} 
\[
\sum_{i,j}(-1)^{i+j}\binom ni\binom nj\binom{2n+j}n\binom{n+j}n\binom
n{2j-i}=\binom{2n}n^3 
\]
\[
\]

\textbf{(Lxxxiv)} 
\[
\sum_{i,j}(-1)^{i+j}\binom ni\binom nj\binom{i+j}n\binom{3n-i-j}n=(-1)^n%
\binom{2n}n 
\]
\[
\]

\textbf{(Lxxxv)} 
\[
\sum_{i,j}(-1)^{i+j}\binom ni\binom nj\binom{i+j}j\binom{n+j}n\binom{2j}j=%
\binom{2n}n^2 
\]
\[
\]

\textbf{(Lxxxvi)} 
\[
\sum_{i,j}(-1)^{i+j}\binom{n+j}n\binom{2n+i}n\binom{2n+j}n\binom{3n+1}{i-j}%
=(-1)^n\binom{2n}n 
\]
\[
\]

\textbf{(Lxxxvii)} 
\[
\sum_{i,j}(-1)^{n+i+j}\binom{2n-i}n^k\binom{n+j}n\binom{n+i}{n-j}=\binom{2n}%
n^k 
\]
\[
\]

\textbf{(Lxxxviii)} 
\[
\sum_{i,j}(-1)^{i+j}\binom ni\binom nj\binom{i+j}n\binom{2i+2j}n=2^n 
\]
\[
\]

\textbf{(Lxxxix)} 
\[
\sum_{i,j}(-1)^{i+j}\binom ni\binom nj\binom{i+j}n^2=\binom{2n}n 
\]
\[
\]

\textbf{(Lxxxx))} 
\[
\sum_{i,j}(-1)^{i+j}\binom{2n-i}n^3\binom{n+j}n\binom{2n+j}n\binom{2n+1}{i-j}%
=\binom{2n}n^4 
\]
\[
\]

(\textbf{Lxxxxi)} 
\[
\sum_{i,j}(-1)^{i+j}\binom{2n-i}n\binom{n+j}n^2\binom{2n+j}n\binom{3n+1}{i-j}%
=\binom{2n}n^3 
\]
\[
\]

\textbf{(Lxxxxi)} 
\[
\sum_{i,j}(-1)^{i+j}\binom{2n-i}n\binom{n+j}n\binom{2n+j}n\binom{2n+1}{i-j}=%
\binom{2n}n^2 
\]
\[
\]

\textbf{Lxxxxii)} 
\[
\sum_{i,j}(-1)^{i+j}\binom{2n-i}n\binom{n+j}n\binom{3n-i}n\binom{2n+j}n%
\binom{3n+1}{i-j}=\binom{2n}n^3 
\]
\[
\]

\textbf{(Lxxxxiii)} 
\[
\sum_{i,j}(-1)^{i+j}\binom{i+j}n\binom{2n-i}n\binom{n+j}n\binom{3n-i}n\binom{%
2n+j}n\binom{3n}{i+j}=\binom{2n}n^2\binom{3n}n^2 
\]
\[
\]

\textbf{(Lxxxxiv)} 
\[
\sum_{i,j}(-1)^{i+j}\binom{2n-i}n\binom{n+j}n\binom{2n+i}n\binom{3n+1}{i-j}%
=(-1)^n\binom{2n}n 
\]
\[
\]

\textbf{(Lxxxxv)} 
\[
\sum_{i,j}(-1)^{i+j}\binom{i+j}n\binom{2n-i}n\binom{n+j}n\binom{3n}{i+j}=-n%
\binom{3n}n 
\]
\[
\]

\textbf{(Lxxxxvi)} 
\[
\sum_{i,j}(-1)^{i+j}\binom{i+j}n\binom{2n-i}n\binom{n+j}n\binom{2n}{i+j}%
=(-1)^n\binom{2n}n\binom{2n+1}n 
\]
\[
\]

\textbf{(Lxxxxvii)} 
\[
\sum_{i,j}(-1)^{i+j}\binom ni\binom nj\binom{2n}{i+j}^2=\binom{2n}n\binom{3n}%
n 
\]
\[
\]

\textbf{(Lxxxxviii)} 
\[
\sum_{i,j}(-1)^{i+j}\binom ni\binom nj\binom{n+i+j}n^2=\binom{2n}n 
\]
\[
\]

\textbf{(Lxxxxviiii)} 
\[
\sum_{i,j}(-1)^{i+j}\binom ni\binom nj\binom{i+j}n^2\binom{n+i+j}n=\binom{2n}%
n^3 
\]
\[
\]

\textbf{(C)} 
\[
\sum_{i,j}\binom ni^2\binom nj\binom ij\binom{3n-j}n=\binom{2n}n^3 
\]
\[
\]

\textbf{(Ci)} 
\[
\sum_{i,j}\binom ni^2\binom nj\binom{n+i-j}{n-j}\binom{2n}{n+j}=\binom{2n}%
n^3 
\]
\[
\]

\textbf{(Cii)} 
\[
\sum_{i,j}\binom ni^2\binom nj\binom{n+i}{n+j}=\binom{2n}n^2 
\]
\[
\]

\textbf{(Ciii)} 
\[
\sum_{i,j}(-1)^{i+j}\binom ni\binom nj\binom ij\binom{2n-j}n\binom{n+j}n%
\binom{3j}{2n}=\binom{2n}n\binom{3n}n 
\]
\[
\]

\textbf{(Civ)} 
\[
\sum_{i,j}(-1)^{i+j}\binom ni\binom nj\binom{n+i}n^3\binom{2j}i=2^n\binom{2n}%
n^3 
\]
\[
\]

\textbf{(Cv)} 
\[
\sum_{i,j}(-1)^{n+i+j}\binom ni^2\binom nj\binom{2j}i=2^n 
\]
\[
\]

\textbf{(Cvi)} 
\[
\sum_{i,j}(-1)^{n+i+j}\binom ni\binom nj\binom{n+j}n^2\binom{n+i-j}{n-j}=1 
\]
\[
\]

\textbf{(Cvii)} 
\[
\sum_{i,j}(-1)^{n+i+j}\binom ni\binom nj\binom{n+j}n^2\binom{n+i}{n-j}=1 
\]
\[
\]

\textbf{(Cviii)} 
\[
\sum_{i,j}(-1)^{n+i+j}\binom ni\binom nj\binom{n+j}n^2\binom{2n+i-j}n=1 
\]
\[
\]

\textbf{(Cix)} 
\[
\sum_{i,j}(-1)^{n+i+j}\binom ni\binom nj\binom{n+j}n\binom{2n-j}n\binom{%
2n+i-j}{n-j}=\binom{2n}n 
\]
\[
\]

\textbf{(Cx)} 
\[
\sum_{i,j}(-1)^{i+j}\binom ni\binom nj\binom{n+j}n\binom{2n-j}n\binom{2n+i}j=%
\binom{2n}n 
\]
\[
\]

\textbf{(Cxi)} 
\[
\sum_k\binom{n+k}n^2\binom{2n+k}n\binom{3n+1}{n-2k}=\binom{2n}n^3 
\]
\[
\]

\textbf{(Cxii)} 
\[
\sum_{i,j}\binom ni\binom nj\binom{i+j}n\binom{3n}{i+j}=\binom{2n}n\binom{4n%
}{2n} 
\]
\[
\]

\textbf{(Cxiii)} 
\[
\sum_{i,j}(-1)^{i+j}\binom ni\binom nj\binom{2n-j}n\binom{3n-j}n\binom
n{2j-i}=\binom{2n}n^3 
\]
\[
\]

\textbf{(Cxiv)} 
\[
\sum_{i,j}(-1)^{i+j}\binom ni\binom nj\binom ij\binom{2n+j}n\binom{3j}{2n}=%
\binom{3n}n^2 
\]
\[
\]

\textbf{(Cxv)} 
\[
\sum_{i,j}(-1)^{i+j}\binom ni\binom nj\binom ij\binom{3n+j}n\binom{3j}{2n}=%
\binom{2n}n\binom{4n}{2n} 
\]
\[
\]

\textbf{(Cxvi)} 
\[
\sum_{i,j}(-1)^{i+j}\binom ni\binom nj\binom{i+j}j\binom{3j}{2n}=\binom{3n}n 
\]
\[
\]

\textbf{(Cxvii)} 
\[
\sum_{i,j}(-1)^{i+j}\binom ni\binom nj\binom{i+j}j\binom{2j}n\binom{3j}{2n}=%
\binom{2n}n\binom{3n}n 
\]
\[
\]

\textbf{(Cxviii)} 
\[
\sum_{i,j}(-1)^{i+j}\binom ni\binom nj^2\binom{n+j}n\binom{3i+1}n=3^n\binom{%
2n}n 
\]
\[
\]

\textbf{(Cxix)} 
\[
\sum_{i,j}(-1)^{i+j}\binom ni\binom nj^2\binom{n+j}n\binom{2n+j}n\binom{3i+1}%
n=3^n\binom{2n}n\binom{3n}n 
\]
\[
\]

\textbf{(Cxx)} 
\[
\sum_{i,j}(-1)^{i+j}\binom ni\binom nj\binom{2n-j}n\binom{3i+1}j=3^n 
\]
\[
\]

\textbf{(Cxxi)} 
\[
\sum_{i,j}(-1)^{i+j}\binom ni\binom nj\binom{n+j}n\binom{n+i+j}n=1 
\]
\[
\]

\textbf{(Cxxii)} 
\[
\sum_{i,j}(-1)^{i+j}\binom{2n-i}n\binom{n+j}n\binom{2n+j}n\binom{3n-i}n%
\binom{4n+1}{i-j}=\binom{2n}n^2 
\]
\[
\]

\textbf{(Cxxiii)} 
\[
\sum_k\binom nk^2\binom{n+k}n\binom{2n+k}n=\binom{2n}n^3 
\]
\[
\]

\textbf{(Cxxiv)} 
\[
\sum_{i,j}(-1)^{i+j}\binom ni\binom nj\binom{i+j}j\binom{n+j}n^k\binom{3j}{2n%
}=\binom{3n}n^{k+1} 
\]
\[
\]

\textbf{(Cxxv)} 
\[
\sum_{i,j}(-1)^{i+j}\binom ni\binom nj\binom ij\binom{3n-j}n\binom{3j}{2n}=%
\binom{2n}n\binom{3n}n 
\]
\[
\]

\textbf{(Cxxvi)} 
\[
\sum_{i,j}(-1)^{i+j}\binom ni\binom nj\binom ij\binom{3n-j}n^2\binom{3j}{2n}=%
\binom{2n}n^2\binom{3n}n 
\]
\[
\]

\textbf{(Cxxvii)} 
\[
\sum_{i,j}(-1)^{i+j}\binom ni\binom nj\binom{i+j}j\binom{2n-j}n^k\binom{3j}{%
2n}=\binom{3n}n 
\]
\[
\]

\textbf{(Cxxviii)} 
\[
\sum_{i,j}\binom ni^2\binom nj\binom{n+i}{n+j}=\binom{2n}n^2 
\]
\[
\]

\textbf{(Cxxix)} 
\[
\sum_{i,j}\binom ni\binom nj\binom{2n+j}n\binom{2j}j\binom{2j}{i-j}=\binom{2n%
}n^3 
\]
\[
\]

\textbf{(Cxxx)} 
\[
\sum_{i,j}(-1)^{i+j}\binom{n+i}i^2\binom{2n-j}n\binom{2n+i}n\binom{3n-j}n%
\binom{4n+1}{j-i}=\binom{2n}n^4 
\]
\[
\]

\textbf{(Cxxxi)} 
\[
\sum_{i,j}(-1)^{i+j}\binom ni\binom nj\binom{3n+j}{2n}\binom n{2j-i}=\binom{%
3n}n^2 
\]
\[
\]

\textbf{(Cxxxii)} 
\[
\sum_{i,j}(-1)^{i+j}\binom ni\binom nj\binom{3n-j}{2n}\binom n{2j-i}=\binom{%
2n}n^2 
\]
\[
\]

\textbf{(Cxxxiii)} 
\[
\sum_{i,j}(-1)^{i+j}\binom ni\binom nj\binom{2n+j}{2n}\binom n{2j-i}=\binom{%
2n}n^2 
\]
\[
\]

\textbf{(Cxxxiv)} 
\[
\sum_{i,j}(-1)^{i+j}\binom{2n-i}n\binom{n+j}j^2\binom{2n+j}n\binom{3n}{i-j}=%
\binom{2n}n^3 
\]

\textbf{(Cxxxv)} 
\[
\sum_{i,j}(-1)^{i+j}\binom{2n-i}n\binom{n+j}j\binom{2n+j}n\binom{3n}{i-j}=%
\binom{2n}n 
\]
\[
\]

\textbf{(Cxxxvi)} 
\[
\sum_{i,j}\binom ni\binom nj\binom{n+j}n^2\binom{n+i}n\binom{2n+i}n\binom
n{i-j}=\binom{2n}n^5 
\]
\[
\]

\textbf{(Cxxxvii)} 
\[
\sum_k\binom nk^4\binom{2n+k}n\binom{3n-k}n\binom{3n+k}n\binom{4n-k}n\binom{%
2k}k^{-1}\binom{2n-2k}{n-k}^{-1}\binom{2n}{2k}^{-1}=\binom{3n}n^2 
\]
\[
\]

The following sums are zero for odd $n$ and the given value for even $n.$

\textbf{(aa)} 
\[
\sum_{i,j}(-1)^{i+j}\binom{2n-i}n\binom{n+j}n\binom{i+j}n\binom{3n+1}{i-j}=%
\binom{3n/2}{n/2} 
\]

\[
\]

\textbf{(bb)} 
\[
\sum_{i,j}(-1)^{i+j}\binom{2n-i}n\binom{2n-j}n\binom{2n}j\binom{2n+1}{i-j}%
=(-1)^{n/2}\binom{2n}n\binom n{n/2} 
\]
\[
\]

\textbf{(cc)} 
\[
\sum_k(-1)^{n/2+k}\binom nk\binom{n+k}n\binom{2n-k}n=\binom{3n/2}{n/2}\binom
n{n/2} 
\]
\[
\]

\textbf{(dd)} 
\[
\sum_k(-1)^{n+k}\binom nk\binom{n+k}n\binom{2n-2k}{n-k}=\binom{2n/3}{n/3}%
\binom{4n/3}{2n/3}\text{ } 
\]
\[
\text{if }3\mid n\text{, }0\text{ otherwise} 
\]

\textbf{(ee)} 
\[
\sum_{i,j}(-1)^{i+j}\binom ni\binom nj\binom{3n-i}n\binom{2n-j}n\binom{n+i-j}%
n=(-1)^{n/2}\binom{2n}n\binom n{n/2}\binom{3n/2}{n/2} 
\]
\[
\]

\textbf{(ff)} 
\[
\sum_{i,j}\binom ni^2\binom nj^2\binom{2i}{n+j}\binom{3n-2i-j}{2n-j}=\binom
n{n/2}^2 
\]
\[
\]

\textbf{(gg)} 
\[
\sum_{i,j}(-1)^{i+j}\binom ni\binom nj\binom{2i}{i-j}=\binom n{n/3} 
\]
\[
\text{if }3\mid n\text{, }0\text{ otherwise} 
\]
\[
\]

\textbf{(hh)} 
\[
\sum_{i,j}(-1)^{i+j}\binom ni^2\binom nj\binom{i+j}n=(-1)^n\binom n{n/2} 
\]
\[
\]

\textbf{(ii)} 
\[
\binom{2n}n\sum_{i,j}(-1)^{n+i}\binom ni\binom ij\binom{2n+i}n\binom{3n+i}n=%
\binom{3n}n\binom{2n}{n/2} 
\]
\[
\]

\textbf{(jj)} 
\[
\sum_{i,j}(-1)^{i+j}\binom ni\binom nj\binom{n+j}n\binom{2n-j}n\binom{i+j}%
n=(-1)^{n/2}\binom n{n/2}\binom{3n/2}{n/2} 
\]
\[
\]

\textbf{(kk)} 
\[
\sum_{i,j}(-1)^j\binom ni^2\binom{i+j}n\binom{n+i+j}n\binom{3n+1}{n-j}=%
\binom{2n}n\binom n{n/2}\binom{3n/2}{n/2} 
\]
\[
\]

\textbf{Triple sums.}

\textbf{18.} 
\[
A_n=\binom{2n}n\sum_{i,j,k}(-1)^{i+j}\binom ni^2\binom nj\binom nk\binom
ij\binom jk\binom n{k-i} 
\]
\[
\]

\textbf{19.} 
\[
A_n=\sum_{i,j,k}\binom ni^2\binom nj^2\binom nk\binom ij\binom jk 
\]
\[
A_n=\sum_{i,j,k}\binom ni^2\binom nj\binom nk\binom ij\binom ik\binom{2j}n 
\]
\[
A_n=\sum_{i,j,k}\binom ni^2\binom nj\binom nk\binom ij\binom ik\binom{2k}n 
\]
\[
\]

\textbf{26.} 
\[
A_n=\binom{2n}n\sum_{i,j,k}\binom ni^2\binom nj\binom nk\binom ij\binom jk 
\]
\[
A_n=\sum_{i,j,k}\binom ni^2\binom nj\binom nk\binom ij\binom jk\binom{2n+j-k%
}{n-k} 
\]
\[
A_n=\sum_{i,j,k}\binom ni^2\binom nj\binom nk\binom ij\binom{j+k}k\binom{n+j%
}{n-k} 
\]
\[
\]

\textbf{27.} 
\[
A_n=\sum_{i,j,k}\binom ni^2\binom nj\binom nk\binom ij\binom jk\binom{2n-i}n 
\]
\[
A_n=\sum_{i,j,k}\binom ni^2\binom nj\binom nk\binom ij\binom jk\binom{n+j-k}%
j 
\]
\[
A_n=\sum_{i,j,k}\binom ni^2\binom nj\binom nk\binom ij\binom jk\binom{n+k}n 
\]
\[
A_n=\sum_{i,j,k}\binom ni^2\binom nj^2\binom nk\binom jk\binom{n+i-j}n 
\]
\[
A_n=\sum_{i,j,k}\binom ni^2\binom nj^2\binom nk\binom ik\binom{j+k}n 
\]
\[
\]

\textbf{28.} 
\[
A_n=\sum_{i,j,k}\binom ni^2\binom nj\binom nk\binom ij\binom jk\binom{n+i-k}%
n 
\]
\[
A_n=\sum_{i,j,k}\binom ni^2\binom nj\binom nk\binom ij\binom jk\binom{i+j}j 
\]
\[
A_n=\sum_{i,j,k}\binom ni^2\binom nj\binom nk\binom ik\binom jk\binom{j+j}k 
\]

\textbf{29.} 
\[
A_n=\binom{2n}n\sum_{i,j,k}\binom ni^2\binom nj\binom nk\binom ij\binom ik 
\]
\[
A_n=\sum_{i,j,k}\binom ni^2\binom nj^2\binom nk\binom ik\binom{n+i}n 
\]
\[
A_n=\sum_{i,j,k}\binom ni^2\binom nj\binom nk\binom ij\binom jk\binom{2n+i-k}%
n 
\]
\[
A_n=\sum_{i,j,k}\binom ni^2\binom nj\binom nk\binom ij\binom jk\binom{2n+i-k%
}{n-k} 
\]
\[
A_n=\sum_{i,j,k}\binom ni^2\binom nj\binom nk\binom ij\binom{i+k}k\binom{n+i%
}{n-k} 
\]
\[
\]

\textbf{101.} 
\[
A_n=\sum_{i,j,k}\binom ni^2\binom nj\binom nk\binom ij\binom{i+k}k\binom{n+j%
}{n-k} 
\]
\[
\]

\textbf{189.} 
\[
A_n^{\prime }=\sum_{i,j,k}\binom ni^2\binom nj\binom nk\binom ij\binom ik%
\binom{j+k}k 
\]
\[
\]

\textbf{193.} 
\[
A_n=\sum_{i,j,k}\binom ni^2\binom nj\binom nk\binom ij\binom{n+k}n\binom{n+i%
}{n-k} 
\]
\[
\]

\textbf{195.} 
\[
A_n=\sum_{i,j,k}\binom ni^2\binom nj\binom nk\binom ij\binom ik\binom{n+j+k}%
j 
\]
\[
\]

\textbf{196.} 
\[
A_n=\sum_{i,j,k}\binom ni^2\binom nj\binom nk\binom ij\binom jk\binom{n+j-k}%
n 
\]
\[
A_n=\sum_{i,j,k}\binom ni^2\binom nj\binom nk\binom ij\binom jk\binom{i+k}k 
\]
\[
\]

\textbf{197.} 
\[
A_n=\sum_{i,j,k}\binom ni\binom nj\binom nk\binom ij\binom ik\binom jk\binom{%
j+k}i 
\]
\[
\]

\textbf{198.} 
\[
A_n=\sum_{i,j,k}\binom ni^2\binom nj\binom nk\binom ij\binom jk\binom{n+k}n 
\]
\[
A_n=\sum_{i,j,k}\binom ni^2\binom nj\binom nk\binom ij\binom jk\binom{2n-i-k}%
n 
\]
\[
A_n=\sum_{i,j,k}\binom ni^2\binom nj\binom nk\binom ij\binom jk\binom{2n-k}n 
\]
\[
A_n=\sum_{i,j,k}\binom ni^2\binom nj^2\binom nk\binom jk\binom{i+k}k 
\]
\[
A_n=\sum_{i,j,k}\binom ni^2\binom nj^2\binom nk\binom ik\binom{i+j}n 
\]
\[
A_n=\sum_{i,j,k}\binom ni^2\binom nj\binom nk\binom ij\binom ik\binom{n+j}n 
\]
\[
\]

\textbf{202.} 
\[
A_n=\sum_{i,j,k}\binom ni^2\binom nj\binom nk\binom ij\binom jk\binom{n+i-k}{%
n-k} 
\]
\[
A_n=\sum_{i,j,k}\binom ni^2\binom nj\binom nk\binom ij\binom ik\binom{n+j-k}{%
n-k} 
\]
\[
A_n=\sum_{i,j,k}\binom ni^2\binom nj\binom nk\binom ij\binom jk\binom{2i}k 
\]
\[
\]

\textbf{212.} 
\[
A_n=\sum_{i,j,k}\binom ni^2\binom nj\binom nk\binom ij\binom jk\binom{i+k}n 
\]
\[
A_n=\sum_{i,j,k}\binom ni^2\binom nj\binom nk\binom ij\binom ik\binom{j+k}n 
\]
\[
A_n=\sum_{i,j,k}\binom ni^2\binom nj\binom nk\binom ij\binom jk\binom{n+i-k}%
i 
\]
\[
A_n=\sum_{i,j,k}\binom ni^2\binom nj\binom nk\binom ij\binom jk\binom{2n+i-k}%
n 
\]
\[
\]

\textbf{Coefficients coming from Laurent polynoms.}

\smallskip Here $c.t.$ means constant term.

\smallskip

\textbf{1.} 
\[
A_n=c.t.\left\{ (x+y+z+t+\frac 1{xyzt})^{5n}\right\} 
\]

\textbf{2.} 
\[
A_n=c.t.\left\{ (x+y+z+t+\frac 1{x^5y^2zt})^{10n}\right\} 
\]

\textbf{3.} 
\[
A_n=c.t.\left\{ 
\begin{array}{c}
(x(1+\frac 1y+\frac 1z+\frac 1t+\frac 1{yz}+\frac 1{yt}+\frac 1{zt}+\frac
1{yzt}) \\ 
+\frac 1x(1+y+z+t+yz+yt+zt+yzt))^{2n}%
\end{array}
\right\} 
\]

\textbf{4.} 
\[
A_n=c.t.\left\{ \frac 1x(1+y+z+t+yt+zt)+\frac{x^2}t(\frac 1y+\frac 1z+\frac
1t))^{3n}\right\} 
\]

\textbf{5.} 
\[
A_n=c.t.\left\{ (\frac 1x(1+y+z+t+yz+zt)+\frac{x^2}{yz}(1+\frac
1t))^{3n}\right\} 
\]

\textbf{6.} 
\[
A_n=c.t.\left\{ (x+y+z+t+\frac 1{xyz}+\frac 1{xyt})^{4n}\right\} 
\]

\textbf{7.} 
\[
A_n=c.t.\left\{ (x+xy^2+xz^8+xt^4+\frac 1{x^7y^8z^8t^4})^{8n}\right\} 
\]

\textbf{8.} 
\[
A_n=c.t.\left\{ (x+y+z+t+\frac 1{x^2yzt})^{6n}\right\} 
\]

\textbf{10.} 
\[
A_n=c.t.\left\{ (\frac 1x+\frac xy+\frac zx+\frac tx+\frac{xt}y+\frac{xz}y+%
\frac{xy}{zt}+\frac{y^2}{xzt})^{4n}\right\} 
\]

\textbf{11.} 
\[
A_n=c.t.\left\{ (\frac 1x+\frac xy+\frac zx+\frac tx+\frac{xt}y+\frac{xz}y+%
\frac{y^2}{xzt})^{4n}\right\} 
\]

\textbf{12.} 
\[
A_n=c.t.\left\{ (x+y+z+t+\frac 1{xyzt}(\frac 1x+\frac 1y))^{6n}\right\} 
\]

\textbf{14.} 
\[
A_n=c.t.\left\{ (x+y+z+t+\frac 1{x^3y}(\frac 1z+\frac 1t))^{6n}\right\} 
\]

\textbf{16} 
\[
A_n=c.t.\left\{ (x+y+z+t+\frac 1x+\frac 1y+\frac 1z+\frac 1t))^{2n}\right\} 
\]

\textbf{24.} 
\[
A_n=c.t.\left\{ (x+y+z+t+\frac 1t(\frac 1x+\frac 1y+\frac 1z)+\frac
1{xyzt^2})^{3n}\right\} 
\]

\textbf{25.} 
\[
A_n=c.t.\left\{ 
\begin{array}{c}
(\frac 1x+\frac 1x(y+z+t)+x(\frac 1y+\frac 1z+\frac 1t) \\ 
\frac 1x(yt+yz+zt)+x(\frac 1{yz}+\frac 1{yt}+\frac 1{zt})+\frac x{yzt})^{2n}%
\end{array}
\right\} 
\]

\textbf{26.} 
\[
A_n=c.t.\left\{ (x+y+z+t+\frac 1x+\frac 1y+\frac 1z+\frac 1t+\frac
y{xt}+\frac z{xt}+\frac t{yz})^{2n}\right\} 
\]

\textbf{29.} 
\[
A_n=c.t.\left\{ 
\begin{array}{c}
(\frac 1x+\frac xy+\frac yx+\frac xz+\frac zx+\frac xt+\frac tx \\ 
+\frac x{yz}+\frac{yz}x+\frac x{yt}+\frac{yt}x+\frac{zt}x)^{2n}%
\end{array}
\right\} 
\]

\textbf{42.} 
\[
A_n=c.t.\left\{ (x+y+z+t+\frac 1x+\frac 1y+\frac 1z+\frac 1t+\frac{zt}%
x+\frac 1{yzt})^{2n}\right\} 
\]

\textbf{51.} 
\[
A_n=c.t.\left\{ (x+y+z+t+\frac 1{t^2}(\frac 1x+\frac 1y+\frac 1z)+\frac
1{xyzt^4})^{3n}\right\} 
\]

\textbf{70.} 
\[
A_n=c.t.\left\{ (x+y+z+t+\frac 1t(\frac 1x+\frac 1y+\frac 1z))^{3n}\right\} 
\]

$\widetilde{\mathbf{101}}.$%
\[
A_n=c.t.\left\{ (\frac 1x+xy+\frac 1{xy}+\frac yz+\frac zy+\frac yt+\frac
ty+\frac y{zt}+\frac{zt}y+\frac z{xy}+\frac t{xy}+\frac
y{xzt}))^{2n}\right\} 
\]

\textbf{185.} 
\[
A_n=c.t.\left\{ (x+y+z+t+\frac 1t(\frac 1{xz}+\frac 1{xy}+\frac
1{yz}))^{2n}\right\} 
\]

\textbf{206.} 
\[
A_n=c.t.\left\{ 
\begin{array}{c}
(x+y+z+t+\frac{zt}y+\frac 1x(\frac yz+\frac zy+\frac yt+\frac ty) \\ 
+\frac 1{x^2}(\frac 1y+\frac 1z+\frac 1t)+\frac 1{x^2zt}(y+\frac 1x))^{2n}%
\end{array}
\right\} 
\]

\textbf{209.} 
\[
A_n=c.t.\left\{ 
\begin{array}{c}
(x(\frac 1y+\frac 1z+\frac 1t+\frac 1{yz}+\frac 1{yt}+\frac 1{zt}+\frac
1{yzt}) \\ 
+\frac 1x(1+y+z+t+yz+yt+zt+yzt))^{2n}%
\end{array}
\right\} 
\]

\textbf{214.} 
\[
A_n=c.t.\left\{ 
\begin{array}{c}
(x(\frac 1y+\frac 1z+\frac 1t+\frac 1{yz}+\frac 1{yt}+\frac 1{zt}) \\ 
+\frac 1x(1+y+z+t+yz+yt+zt+yzt))^{2n}%
\end{array}
\right\} 
\]

\textbf{218.} 
\[
A_n=c.t.\left\{ (\frac 1{xy}+\frac xz+\frac zx+\frac xt+\frac tx+\frac{xy}%
z+\frac z{xy}+\frac{xy}t+\frac t{xy}+\frac{zt}x+\frac x{zt}+\frac x{yzt}+%
\frac{yzt}x)^{2n}\right\} 
\]

\textbf{287.} 
\[
A_n=c.t.\left\{ (\frac 1x+\frac yx+\frac xz+\frac zx+\frac xt+\frac tx+\frac
x{yz}+\frac{yz}x+\frac x{yt}+\frac{yt}x+\frac{zt}x+\frac x{yzt}+\frac{yzt}%
x)^{2n}\right\} 
\]

\textbf{308} 
\[
A_n=c.t.\left\{ (x+y+z+t+\frac 1x+\frac 1y+\frac 1z+\frac x{yz}+\frac
y{xt}+\frac z{xt}))^{2n}\right\} 
\]

\textbf{309.} 
\[
A_n=c.t.\left\{ 
\begin{array}{c}
(x(\frac 1y+\frac 1z+\frac 1t+\frac 1{yz}+\frac 1{yt}+\frac 1{yzt}) \\ 
+\frac 1x(1+y+z+t+yz+yt+zt+yzt))^{2n}%
\end{array}
\right\} 
\]
\[
\]

\textbf{References. }

1.G.Almkvist, C.van Enckevort, D.van Straten and W.Zudilin, \textsl{Tables
of Calabi-Yau equations,} math.AG/0507430, (2005).

2.G.Almkvist, C.Krattenthaler,\textsl{\ Some harmonic sums related to
Calabi-Yau differential equations, in preparation. }

3.G.Almkvist and W.Zudilin, \textsl{Differential equations, mirror maps and
zeta values, }math.AG/0507430

4.W.Zudilin, \textsl{Approximation to -, di- and trilogarithms, }%
math.CA/0409023.

\bigskip 

Math. Dept.

Univ. of Lund

Box 118 

22100 Lund, Sweden

gert@maths.lth.se

\end{document}